\documentclass[11pt,a4paper]{amsart}

\usepackage{amsmath}
\usepackage{amssymb}
\usepackage{latexsym}
\usepackage{xypic}

\newcommand{\dual}{\makebox[0mm]{}^{{\scriptstyle\vee}}}
\newtheorem{theorem}{Theorem}[section]
\newtheorem{lemma}[theorem]{Lemma}
\newtheorem{proposition}[theorem]{Proposition}
\newtheorem{definition}[theorem]{Definition}
\newtheorem{corollary}[theorem]{Corollary}
\newtheorem{conjecture}[theorem]{Conjecture}
\newtheorem{exmp}[theorem]{Example}
\newtheorem{exmps}[theorem]{Examples}
\newtheorem{rem}[theorem]{Remark}
\newenvironment{example}{\begin{exmp}\rm}{\end{exmp}}

\newenvironment{remark}{\begin{rem}\rm}{\end{rem}\rm}

\newcommand{\qqed}{\hspace*{\fill}$\Box$}

\newcommand{\beeq}[1]{\begin{eqnarray}\label{#1}}
\newcommand{\eneq}{\end{eqnarray}}
\newcommand{\cal}{\mathcal}

\newcommand{\kc}{{\cal C}}

\newcommand{\ke}{{\cal E}}
\newcommand{\kf}{{\cal F}}

\newcommand{\kl}{{\cal L}}

\newcommand{\ko}{{\cal O}}

\newcommand{\Db}{{\rm D}^{\rm b}}

\newcommand{\OO}{{\rm O}}

\newcommand{\IC}{{\mathbb C}}

\newcommand{\IP}{{\mathbb P}}
\newcommand{\IQ}{{\mathbb Q}}
\newcommand{\IR}{{\mathbb R}}
\newcommand{\IZ}{{\mathbb Z}}

\newcommand{\Hom}{{\rm Hom}}
\newcommand{\Ext}{{\rm Ext}}

\newcommand{\rk}{{\rm rk}}

\newcommand{\verylongarrow}[1]{\hbox to #1{\rightarrowfill}}

\newcommand\mynote[1]

\begin{document}

\title[Equivalences of twisted K3 surfaces]{Equivalences of twisted K3 surfaces}

\author[Daniel Huybrechts]{Daniel Huybrechts}
\address{D.H.: Mathematisches Institut, Universit{\"a}t Bonn,
  Beringstr.\ 1, D- 53115 Bonn, Germany}
\email{huybrech@math.uni-bonn.de}

\author[Paolo Stellari]{Paolo Stellari}
\address{P.S.: Department of Mathematics, Universit{\`a} degli Studi di Milano, Via Cesare
Saldini 50, 20133 Milan, Italy} \email{stellari@mat.unimi.it}

\begin{abstract}
We  prove that two derived equivalent twisted K3 surfaces have
isomorphic periods. The converse is shown for K3 surfaces with
large Picard number. It is also shown that all possible twisted
derived equivalences between arbitrary twisted K3 surfaces form a
subgroup of the group of all orthogonal transformations of the
cohomology of a K3 surface.

The passage from twisted derived equivalences to an action on the
cohomology is made possible by twisted Chern characters that will
be introduced for arbitrary smooth projective varieties.
\end{abstract}

\maketitle

By definition a K3 surface is a compact complex surface $X$ with
trivial  canonical bundle and vanishing $H^1(X,\ko_X)$. As was
shown by Kodaira in \cite{Kod} all K3 surfaces are deformation
equivalent. In particular, any K3 surface is diffeomorphic to the
four-dimensional manifold $M$ underlying the Fermat quartic in
$\IP^3$ defined by  $x_0^4+x_1^4+x_2^4+x_3^4=0$. Thus, we may
think of a K3 surface $X$ as a complex structure $I$ on $M$. (As
it turns out, every complex structure on $M$ does indeed define a
K3 surface, see \cite{FM}.)

In the following, we shall fix the orientation on $M$ that is
induced by a complex structure and denote by $\Lambda$ the
cohomology $H^2(M,\IZ)$ endowed with the intersection pairing.
This is an even unimodular lattice of signature $(3,19)$ and
hence isomorphic to $(-E_8)^{\oplus 2}\oplus U^{\oplus 3}$ with
$E_8$ the unique even positive definite unimodular lattice of rank
eight and $U$ the hyperbolic plane.

We shall denote by $\widetilde \Lambda$ the lattice given by  the
full integral cohomology $H^*(M,\IZ)$ (which is concentrated in
even degree) endowed with the Mukai pairing $\langle
\varphi_0+\varphi_2+\varphi_4,\psi_0+\psi_2+\psi_4\rangle=\varphi_2\wedge\psi_2-
\varphi_0\wedge\psi_4-\varphi_4\wedge\psi_0$. In other words,
$\widetilde \Lambda$ is the direct sum of $(H^0\oplus
H^4)(M,\IZ)$ endowed with the negative intersection pairing and
$\Lambda$. Hence, $\widetilde\Lambda\cong\Lambda\oplus U$.

An \emph{isomorphism} between two K3 surfaces $X$ and $X'$ given
by two complex structures $I$ respectively $I'$ on $M$ is a
diffeomorphism $f\in{\rm Diff}(M)$ such that $I=f^*(I')$. Any
such diffeomorphism $f$ acts on the cohomology of $M$ and,
therefore, induces a lattice automorphism
$f_*:\Lambda\cong\Lambda$.

Conversely, one might wonder whether any element
$g\in\OO(\Lambda)$ is of this form. This is essentially true and
has been proved by Borcea in \cite{Borcea}. The precise statement
is:
\smallskip

{\it For any $g\in\OO(\Lambda)$ there exist two K3 surfaces
$X=(M,I)$ and $X'=(M,I')$ and an isomorphism $f:X\cong X'$ with
$f_*=\pm g$.}
\smallskip

(In fact, we can even prescribe the K3 surface $X=(M,I)$, but stated
like this the result compares nicely with Theorem \ref{firstthmIntro}.)
The proof of this fact uses the full theory of K3 surfaces, i.e.\
Global Torelli theorem, surjectivity of the period map, etc.
Donaldson showed in \cite{Donaldson} that the image of ${\rm
Diff}(M)\to\OO(\Lambda)$ is in fact the index two subgroup
$\OO_+(\Lambda)$ of orthogonal transformation preserving the
orientation of the positive directions. (Note that $\OO(\Lambda)$
is generated by $\OO_+(\Lambda)$ and $\pm{\rm id}$.)

In a next step, we consider a more flexible notion of
isomorphisms of K3 surfaces: One says that two K3 surfaces $X$ and
$X'$ are \emph{derived equivalent} if there exists a
Fourier--Mukai equivalence $\Phi:\Db(X)\cong\Db(X')$. Here,
$\Db(X)$ is the bounded derived category of the abelian category
${\bf Coh}(X)$ of coherent sheaves on $X$. (Usually, derived
equivalence is only considered for algebraic K3 surfaces.)

Clearly, any isomorphism between $X$ and $X'$ given by $f\in{\rm
Diff}(M)$ induces a Fourier--Mukai equivalence $\Phi:=Rf_*$. By
results of Mukai and Orlov one knows how to associate to any
Fourier--Mukai equivalence $\Phi:\Db(X)\cong\Db(X')$ an
isomorphism $\Phi_*:H^*(X,\IZ)\cong H^*(X',\IZ)$  or, thinking of
$X$ and $X'$ as complex structures $I$ respectively $I'$ on $M$,
an orthogonal transformation $\Phi_*\in\OO(\widetilde \Lambda)$.
For $\Phi=Rf_*$ as above, this gives back the standard action of
the diffeomorphism $f$ on $H^*(M,\IZ)$. However, for more general
derived equivalences $\Phi$ the induced
$\Phi_*\in\OO(\widetilde\Lambda)$ does not respect the
decomposition $\widetilde\Lambda=\Lambda\oplus U$. So it makes
perfect sense to generalize the above question on the
cohomological action of isomorphisms between K3 surfaces to the
derived setting:

\medskip

{\it Is any element in $\OO_+(\widetilde\Lambda)$ of the form
$\Phi_*$ for some derived equivalence $\Phi$? Do elements of the
form $\Phi_*$ form a subgroup of $\OO_+(\widetilde \Lambda)$?}

\medskip

As will be shown  in Section \ref{TwistedisometriesSection}, the
answer to both questions is negative. So generalizing
isomorphisms of K3 surfaces to derived equivalence seems not very
natural from the cohomological point of view. In fact, one has to
go one step further in order to get a nice cohomological
behaviour.

Instead of allowing only equivalences $\Phi:\Db(X)\cong\Db(X')$
one considers more generally Fourier--Mukai equivalences
$\Phi:\Db(X,\alpha)\cong\Db(X',\alpha')$ of twisted derived
categories. Here, $\Db(X,\alpha)$ is the bounded derived category
of $\alpha$-twisted coherent sheaves on $X$, where $\alpha$ is a
Brauer class, i.e.\ a torsion class in $H^2(X,\ko_X^*)$. (More
details are recalled in Section \ref{TwistedChernSection}.)

By means of twisted Chern characters and twisted Mukai vector, that
will be introduced in full generality for arbitrary smooth
projective varieties in Section \ref{TwistedChernSection}, we will
show how to associate to any twisted Fourier--Mukai equivalence
$\Phi$ an isometry $\Phi_*\in\OO(\widetilde\Lambda)$ (see also
Theorem \ref{theom2}, i)). In fact, $\Phi_*$ depends on the
additional choices of B-field lifts of $\alpha$ and $\alpha'$. This
will be spelled out in detail in Sections \ref{TwistedChernSection}
and \ref{GCYSection}, but for the relation between B-fields and
Brauer classes see the discussion further below.

In the twisted context, the above question has an affirmative
answer (see Proposition \ref{twistedisoformgroup}):

\begin{theorem}\label{firstthmIntro}
For any $g\in\OO_+(\widetilde\Lambda)$ there exists a twisted
Fourier--Mukai equivalence
$\Phi:\Db(X,\alpha)\cong\Db(X',\alpha')$ of algebraic K3 surfaces
such that $\Phi_*=g$.
\end{theorem}

Note that in contrast to the above result of Donaldson,
the K3 surface $X$ has to be chosen carefully here.
The choice of B-field lifts of $\alpha$ and $\alpha'$ are
tacitly assumed. A B-field is by definition a class
in $H^2(X,\IQ)$.
As it will be shown, it suffices actually to consider derived
equivalences between untwisted derived categories, but with a
non-trivial B-field $B'$ turned on. So, studying twisted derived
categories teaches us that even in the untwisted situation a
non-trivial lift, e.g.\ $0\ne B\in H^2(X,\IZ)$, might be
necessary to obtain a clean picture.

Note that at least conjecturally,
$\Phi_*\in\OO_+(\widetilde\Lambda)$ for any Fourier--Mukai
equivalence $\Phi$.

\medskip

It will turn out useful to look at these purely algebraic
questions from a more differential geometric angle. For this, one
uses the generalization of the notion of K3 surfaces, i.e.\
complex structures on $M$, provided by generalized Calabi--Yau
structures on $M$. The relation between twisted derived
equivalence introduced above and generalized Calabi--Yau
structures works perfectly well on the level of cohomology and is
used to formulate our results on twisted derived equivalences. A
deeper understanding of the interplay between these quite
different mathematical structures needs still to be developed.

A generalized Calabi--Yau structure on $M$ consists of an even
closed complex form $\varphi$ such that
$\varphi_2^2-2\varphi_0\wedge\varphi_4$ is a zero four-form and
$\varphi_2\wedge\bar\varphi_2-\varphi_0\wedge\bar\varphi_4-\bar\varphi_0\wedge\varphi_4$
is a volume form. This notion was introduced by Hitchin in
\cite{Hit}, further discussed in Gualtieri's theses
\cite{Gual} and in the case of K3 surfaces in
\cite{HuyGen}.

There are two types of generalized Calabi--Yau structures on $M$:
Either $\varphi=\exp(B)\cdot\sigma=\sigma+B\wedge\sigma$, where
$\sigma$ is a holomorphic volume form on a K3 surface $X=(M,I)$,
or $\varphi=\exp(B+i\omega)$ with $\omega$ a symplectic structure
on $M$. In both cases $B$ is a real closed two-form.

To any generalized Calabi--Yau structure $\varphi$ one naturally
associates a weight two-Hodge structure on $H^*(M,\IZ)$ be
declaring $[\varphi]\IC$ to be the $(2,0)$-part of it. The
$(1,1)$-part is then given as the orthogonal complement of it with
respect to the Mukai pairing. If $\varphi=\exp(B)\cdot\sigma$ with
$\sigma$ a holomorphic volume form on $X=(M,I)$ then we write
$\widetilde H(X,B,\IZ)$ for this Hodge structure. For details see
\cite{HuyGen} or Section \ref{GCYSection}.

Twisted derived categories and generalized Calabi--Yau structures
of the form $\exp(B)\sigma$ are related by the following
construction: If $X=(M,I)$ is a K3 surface, then the cohomology
class of $B$ has its $(0,2)$-part in $H^2(X,\ko_X)$. Via the
exponential map $\exp:\ko_X\to\ko_X^*$ this yields an element
$\alpha_B\in H^2(X,\ko_X^*)$. If $\sigma$ is a holomorphic volume
form on $X$, then $\alpha_B$ only depends on the cohomology class
of  $B\wedge\sigma$. In other words, $\alpha_B$ is naturally
associated to the cohomology class of $\varphi=\exp(B)\cdot\sigma$.
Conversely, as $H^3(M,\IZ)=0$, any class $\alpha\in
H^2(X,\ko_X^*)$ is of the form $\alpha_B$ for some $B$. Moreover,
$\alpha$ is a torsion class if and only if $\alpha=\alpha_B$ for
some $B\in H^2(M,\IQ)$.

\medskip

There are two problems naturally arising in this discussion:

I) {\it Let $X=(M,I)$ be an algebraic K3 surface with a Brauer
class $\alpha$. Describe the image of the following three
homomorphisms

{\rm i)} ${\rm Aut}(X)\to\OO(\Lambda)$,~~ {\rm ii)} ${\rm
Aut}(\Db(X))\to\OO(\widetilde\Lambda)$~{\rm  and} {\rm iii)} ${\rm
Aut}(\Db(X,\alpha))\to\OO(\widetilde\Lambda)$.}

\medskip

As a consequence of the Global Torelli theorem, one can describe
the image of ${\rm Aut}(X)\to\OO(\Lambda)$ as a certain subgroup
of the group of Hodge isometries of the Hodge structure
$H^2(X,\IZ)$. A complete answer to ii) is not yet known, but
using results of Mukai and Orlov it was observed in
\cite{HLOY1,Ploog} that the image of ${\rm
Aut}(\Db(X))\to\OO(\widetilde\Lambda)$ is a subgroup of index at
most  two inside the subgroup of all Hodge isometries of the
lattice $\widetilde H(X,\IZ)$. In fact, it can be shown that the image
contains the subgroup of all Hodge isometries that preserve the
natural orientation of the four positive directions (see Section
\ref{OrPreservSect} for more details).
Already in \cite{Sz}
Szendr\H{o}i argues that every derived equivalence
should preserve the orientation of the positive
directions. In particular, one indeed expects that
the image of  ${\rm
  Aut}(\Db(X))\to\OO(\widetilde\Lambda)$ is the group
of all orientation preserving Hodge isometries.

The last problem iii) is related to
C\u{a}ld\u{a}raru's conjecture, see Remark \ref{CalConjtwistauto}.

\medskip

II) {\it Let $X=(M,I)$ and $X'=(M,I')$ be two algebraic K3
surfaces endowed with Brauer classes $\alpha$ respectively
$\alpha'$. Find cohomological criteria which determine when

{\rm i)} $X\cong X'$, {\rm ii)} $\Db(X)\cong\Db(X')$, or {\rm
iii)} $\Db(X,\alpha)\cong\Db(X',\alpha')$.}

\medskip

The answer to i) is provided by the Global Torelli theorem:
$X\cong X'$ if and only if there exists a Hodge isometry
$H^2(X,\IZ)\cong H^2(X',\IZ)$ (see \cite{Per}). (But not any
Hodge isometry can be lifted to an isomorphism.)

The derived version of it yields an answer to ii): Let $X$ and
$X'$ be two algebraic K3 surfaces. Then $\Db(X)\cong\Db(X')$ if
and only if $\widetilde H(X,\IZ)\cong \widetilde H(X',\IZ)$. This
result is due to Orlov \cite{Or} and relies on techniques
introduced by Mukai \cite{Mu}.

The answer to iii) is supposed to be provided by the following
conjecture formulated, though in a slightly different form, by
C\u{a}ld\u{a}raru in his thesis \cite{Cal}.

\begin{conjecture}
Let $X$ and $X'$ be two algebraic K3 surfaces with rational
B-fields $B$ respectively $B'$ inducing Brauer classes $\alpha$
respectively $\alpha'$. Then there exists a Fourier--Mukai
equivalence $\Db(X,\alpha)\cong\Db(X',\alpha')$ if and only if
there exists a Hodge isometry $\widetilde H(X,B,\IZ)\cong
\widetilde H(X',B',\IZ)$ that respects the natural orientation of
the four positive directions.
\end{conjecture}

\begin{remark}\label{CalConjtwistauto}
A slightly refined version of this conjectures predicts that
actually any Hodge isometry in $\OO_+$ can  be lifted to a twisted
derived equivalence. This would in particular answer I, iii).

C\u{a}ld\u{a}raru actually conjectured that
$\Db(X,\alpha)\cong\Db(X',\alpha')$ if and only if the
transcendental lattices of the twisted Hodge structures are Hodge
isometric. In fact, in the untwisted situation it is easy to see
that any Hodge isometry $T(X)\cong T(X')$ of the transcendental
lattices lifts to a Hodge isometry $\widetilde
H(X,\IZ)\cong\widetilde H(X',\IZ)$. This does not hold any longer
in the twisted case. Thus, as will be explained in detail in
Section \ref{TwistedFMonCohSection}, one has to modify the
original conjecture of C\u{a}ld\u{a}raru's and use the full
twisted Hodge structure $\widetilde H(X,B,\IZ)$ instead of just
its transcendental lattice.
\end{remark}

Using the twisted Chern character introduced in Section
\ref{TwistedChernSection} one can at least prove parts of this
conjecture (cf.\ Propositions \ref{HyieldsD} and \ref{MainSect6}).

\begin{theorem}\label{theom2}
Let $X,B,\alpha,X',B',\alpha'$ as before.

{\rm i)} If $\Phi:\Db(X,\alpha)\cong\Db(X',\alpha')$ is a
Fourier--Mukai equivalence, then there exists a naturally defined
Hodge isometry $\Phi_*^{B,B'}:\widetilde H(X,B,\IZ)\cong
\widetilde H(X',B',\IZ)$.

{\rm ii)}  If the Picard number $\rho(X)$ satisfies
$\rho(X)\geq12$, then for any orientation preserving
Hodge iso\-metry $g:\widetilde
H(X,B,\IZ)\cong \widetilde H(X',B',\IZ)$ there exists a
Fourier--Mukai equi\-va\-lence
$\Phi:\Db(X,\alpha)\cong\Db(X',\alpha_{B'})$ such that
$\Phi_*=g$.
\end{theorem}

The first part will be shown (cf.\ Corollary
\ref{finiteFMpartners}) to imply
\begin{corollary}
Any twisted algebraic K3 surface $(X,\alpha)$ admits only
finitely many Fourier--Mukai partners, i.e.\ there exists only a
finite number of isomorphism classes of twisted K3 surfaces
$(X',\alpha')$ such that one can find a Fourier--Mukai equivalence
$\Db(X,\alpha)\cong\Db(X',\alpha')$.
\end{corollary}

The untwisted case of this corollary had been proved in \cite{BM}.

\begin{remark}
Orlov has proved that in fact any equivalence between untwisted
derived categories is of Fourier--Mukai type. So in the untwisted
case one could just consider equivalences of derived categories.
An analogous result is expected (and in fact proved in Section
\ref{SectionlargePicard} for large Picard number) also in the
twisted situation, but for the time being we have to restrict to
the geometrically relevant case of Fourier--Mukai equivalences.
\end{remark}

\smallskip

Here is an outline of the paper. In Section
\ref{TwistedChernSection} we introduce twisted Chern characters on
arbitrary smooth projective varieties as an additive map from the
K-group of twisted coherent sheaves to rational cohomology and prove
a few basic facts about them. In particular, we will see that the
standard Hodge conjecture in the twisted context is equivalent to
the standard Hodge conjecture.

Section \ref{GCYSection} explains the relation between
generalized Calabi--Yau structures and twisted K3 surfaces, i.e.\
K3 surfaces endowed with a Brauer class.

In Section \ref{BrauerSect} we will study the Brauer group of a
K3 surface by introducing three equivalence relations on it. The
main result asserts the finiteness of each equivalence class
modulo the action of the group of automorphisms.

The general framework associating a natural map on the cohomology
to a twisted Fourier--Mukai transform is
explained in Section \ref{TwistedFMonCohSection}. For K3 surfaces
we will find that a Fourier--Mukai equivalence yields a Hodge
iso\-metry between Hodge structures that are defined by means of
related genera\-lized Calabi--Yau structures. This will in
particular lead to the finiteness result for twisted
Fourier--Mukai partners. The section also contains a detailed
discussion of C\u{a}ld\u{a}raru's conjecture and explains why and
how it has to be modified.

Section \ref{TwistedisometriesSection} contains the proof of
Theorem \ref{firstthmIntro}. This part is based on purely lattice
theoretical considerations and a proof of C\u{a}ld\u{a}raru's
conjecture for large Picard number given in Section
\ref{SectionlargePicard}.

The last section shows that despite the finiteness of
Fourier--Mukai partners, one always finds arbitrarily many twisted
Fourier--Mukai partners and, unlike the untwisted case, the Picard
group may even be chosen large. The examples in this section
illustrate the difference between the twisted and the untwisted
world alluded to in the preceding sections.

\smallskip

{\bf Acknowledgements:} We wish to thank Tony Pantev for comments on
earlier versions of the paper. We are most grateful to Manfred Lehn
and Christian Meyer for their help with Example \ref{counter}.
Thanks also to K.\ Yoshioka who suggested a simplification of our
proof of Proposition \ref{proporpreserv}. During the preparation of
this paper the second named author was partially supported by the
MIUR of the Italian government in the framework of the National
Research Project ``Algebraic Varieties'' (Cofin 2002).

\smallskip

{\bf Note added in Proof:} Using results of Yoshioka \cite{Y} on the existence
and non-triviality of moduli spaces of twisted sheaves on K3 surfaces
we can now prove C\u{a}ld\u{a}raru's conjecture, see \cite{HS2}.
\section{Twisted Chern characters}\label{TwistedChernSection}

In the following, we let $X$ be a smooth projective variety over
$\IC$ and $\alpha\in H^2(X,\ko_X^*)$ be a torsion class, i.e.\ an
element in the Brauer group ${\rm Br}(X)$. The pair $(X,\alpha)$
will sometimes be called a twisted variety (see Definition
\ref{DeftwistedK3}). We may represent $\alpha$ by a \v{C}ech
2-cocycle $\{\alpha_{ijk}\in\Gamma(U_i\cap U_j\cap
U_k,\mathcal{O}^*_X)\}$ with $X=\bigcup_{i\in I} U_i$  an
appropriate  open analytic cover. An \emph{$\alpha$-twisted
(coherent) sheaf} $E$ consists of pairs $(\{E_i\}_{i\in
I},\{\varphi_{ij}\}_{i,j\in I})$ such that the ${E}_i$ are
(coherent) sheaves on $U_i$ and $\varphi_{ij}:{E}_i|_{U_i\cap
U_j}\rightarrow{E}_j|_{U_i\cap U_j}$ are isomorphisms satisfying
the following conditions:

{\rm i)} $\varphi_{ii}=\mathrm{id}$,

{\rm ii)} $\varphi_{ji}=\varphi_{ij}^{-1}$, and

{\rm iii)}
$\varphi_{ij}\circ\varphi_{jk}\circ\varphi_{ki}=\alpha_{ijk}\cdot
\mathrm{id}$.

\begin{definition}
By ${\bf Coh}(X,\alpha)$ we denote the abelian category of
$\alpha$-twisted coherent sheaves. Its K-group is denoted $K(X,\alpha)$.
\end{definition}

A priori, the above definition of ${\bf Coh}(X,\alpha)$ depends on
the chosen \v{C}ech representative of $\alpha$. However, it is not
difficult to see that for two different choices the two abelian
categories are equivalent (see \cite{Cal}).
Note however that the equivalence depends on the additional
choice of $\{\beta_{ij}\in\ko^*(U_i\cap U_j)\}$ satisfying
$\alpha'_{ijk}\cdot\alpha_{ijk}^{-1}=\beta_{ij}\cdot\beta_{jk}\cdot\beta_{ki}$,
where $\{\alpha_{ijk}\}$ and $\{\alpha'_{ijk}\}$ are two \v{C}ech
cocycles
representing the same Brauer class.
Let us also introduce
the notation $\Db(X,\alpha)$ for the bounded derived category of
${\bf Coh}(X,\alpha)$, although we won't say anything about it in
this section.

As the tensor product $\kf\otimes\ke$ of an $\alpha$-twisted sheaf
$\kf$ with a $\beta$-twisted sheaf $\ke$ is an
$\alpha\cdot\beta$-twisted sheaf, the abelian category ${\bf
Coh}(X,\alpha)$ has no natural tensor structure (except when
$\alpha$ is trivial). Hence, its K-group $K(X,\alpha)$ is indeed
just an additive group and not a ring.

The aim of this section is to construct a twisted Chern character
$$\xymatrix{K(X,\alpha)\ar[r]&H^*(X,\IQ)}$$
with a number of specific properties which are all twisted
versions of the standard results on Chern characters. As it will
turn out, however, the definition depends on the additional
choice of a B-field, i.e.\ a cohomology
class in $H^2(X,\IQ)$, and only works for topologically trivial
Brauer classes, i.e.\ for those $\alpha\in H^2(X,\ko_X^*)$ with
trivial boundary in $H^3(X,\IZ)$ under the exponential sequence
$$\xymatrix{0\ar[r]&\IZ\ar[r]&\ko_X\ar[r]^\exp&\ko_X^*\ar[r]&0.}$$
Thus, our construction works for arbitrary Brauer classes on
smooth projective varieties $X$ with $H^3(X,\IZ)_{\rm tor}=0$,
e.g.\ for  K3 surfaces.

\begin{proposition}\label{ProptwistedChern}
Suppose $B\in H^2(X,\IQ)$ is a rational B-field such that its
$(0,2)$-part $B^{0,2}\in H^2(X,\ko_X)$ maps to $\alpha$, i.e.\
$\exp(B^{0,2})=\alpha$.

Then there exists a map
$$\xymatrix{{\rm ch}^B:K(X,\alpha)\ar[r]&H^*(X,\IQ)}$$
such that:

{\rm i)} ${\rm ch}^B$ is additive, i.e.\ ${\rm ch}^B(E_1\oplus
E_2)={\rm ch}^B(E_1)+{\rm ch}^B(E_2)$.

{\rm ii)} If $B={\rm c}_1(L)\in H^2(X,\IZ)$, then ${\rm
ch}^B(E)=\exp({\rm c}_1(L))\cdot {\rm ch}(E)$. (Note that with
this assumption $\alpha$ is trivial and an $\alpha$-twisted sheaf
is just an ordinary sheaf.)

{\rm iii)} For two choices $(B_1,\alpha_1:=\exp(B_1^{0,2}))$,
$(B_2,\alpha_2:=\exp(B_2^{0,2}))$ and $E_i\in K(X,\alpha_i)$ one
has $${\rm ch}^{B_1}(E_1)\cdot {\rm ch}^{B_2}(E_2)={\rm
ch}^{B_1+B_2}(E_1\otimes E_2).$$

{\rm iv)} For any $E\in K(X,\alpha)$ one has ${\rm ch}^B(E)\in
\exp(B)\left(\bigoplus H^{p,p}(X)\right)$.
\end{proposition}

\begin{proof} For the basic facts on twisted sheaves we refer to
\cite{Cal}.

Any coherent $\alpha$-twisted sheaf admits a finite resolution by
locally free sheaves. (Here, one definitely needs  $\alpha$ be
torsion. For the argument see \cite[Lemma 2.1.4]{Cal}.)
Hence, $K(X,\alpha)$ can also be regarded as the
Grothendieck group of locally free $\alpha$-twisted coherent
sheaves. Therefore, it suffices to define ${\rm ch}^B$ for
locally free $\alpha$-twisted sheaves in a way that it becomes
additive for short exact sequences.

Let us fix a \v{C}ech representative
$\{\alpha_{ijk}\in\Gamma(U_{ijk},\ko^*)\}$ of $\alpha\in
H^2(X,\ko^*)$.


In fact, in our situation we may start with a \v{C}ech
representative $B_{ijk}\in\Gamma(U_{ijk},\IQ)$ of the B-field $B$
and use the cocycle given by $\alpha_{ijk}:=\exp(B_{ijk})$
regarded as local sections of $\IR/\IZ={\rm U}(1)\subset\ko^*$
that represents $\alpha$.

Since $\kc^\infty$ is acyclic, there exist functions
$a_{ij}\in\Gamma(U_{ij},\kc^\infty)$ with
$-a_{ij}+a_{ik}-a_{jk}=B_{ijk}$. (We assume that the cover is
sufficiently fine.)

Let now $E$ be an $\alpha$-twisted sheaf (it does not need to be
coherent) given by $\{E_i,\varphi_{ij}\}$. Then consider
$\{E_i,\varphi_{ij}':=\varphi_{ij}\cdot \exp(a_{ij})\}$. It is
easy to check that $\{\varphi_{ij}'\}$ is an honest cocycle, i.e.\
$\varphi_{ij}'\circ\varphi_{jk}'\circ\varphi_{ki}'={\rm id}$.
Hence, we have constructed a(n untwisted) sheaf
$E_B:=\{E_i,\varphi_{ij}'\}$.

Then define
$${\rm ch}^B(E):={\rm ch}(E_B).$$

We have to check that with this definition  ${\rm ch}^B(E)$ does
not depend on any of the choices. This is straightforward and
left to the reader.

i) is also an immediate consequence of the construction and ii)
follows from $E_B=E\otimes L$ under the assumptions. For iii) one
observes that ${E_1}_{B_1}\otimes {E_2}_{B_2}\cong {(E_1\otimes
E_2)}_{B_1+B_2}$.

Suppose $B_0:=k\cdot B\in H^2(X,\IZ)$ for some non-trivial $k\in
\IZ$. Due to ii) and iii) one has for any $E\in K(X,\alpha)$
$${\rm ch}^B(E)^k={\rm ch}^{B_0}(E^{\otimes k})=\exp(B_0)\cdot{\rm
ch}(E^{\otimes k}).$$ Hence,
\begin{eqnarray*} \left(\exp(-B)\cdot{\rm
ch}^B(E)\right)^k&=&\exp(-B_0)\cdot{\rm
ch}^B(E)^k\\
&=&{\rm ch}(E^{\otimes k})\in \bigoplus H^{p,p}(X),
\end{eqnarray*}
for $E^{\otimes k}$ is an algebraic vector bundle and has thus
Chern classes of pure type.

The assertion iv) now follows from the following easy observation:

If $(v_0,v_1,\ldots,v_n)\in H^0\oplus H^{2}\oplus\ldots\oplus
H^{2n}$ with $v_0\ne 0$ and such that $v^k\in\bigoplus H^{p,p}$,
then also $v\in\bigoplus H^{p,p}$, i.e.\ $v_p\in H^{p,p}$ for all
$p$.

Suppose we have shown that  $v_i\in H^{i,i}$ for all $i<j$. Write
$(v^k)_j=k \cdot(v_j v_0)+P(v_0,\ldots,v_{j-1})$ with $P$ a
certain polynomial. By assumption $(v^k)_j$ is of pure type and
by induction hypotheses the same holds for
$P(v_0,\ldots,v_{j-1})$. Since $v_0\ne0$, this yields $v_j\in
H^{j,j}$.
\end{proof}

\begin{remark}\label{MHSRem}{\rm(Mixed Hodge structures by twisting)}

i) The cohomology $H^*(X,\IZ)$ can be seen as a direct sum of
Hodge structures (over $\IZ$) and, therefore, as a mixed Hodge
structure with ascending filtration $W_i$ given by
$W_i=\bigoplus_{j\leq i} H^j(X,\IQ)$. For any $B\in H^2(X,\IQ)$
the isomorphism
$$\exp(B):H^*(X,\IQ)\cong H^*(X,\IQ)$$ induces a natural mixed
Hodge structure given by
$$W_i:=\exp(B)\left(\bigoplus_{j\leq i} H^j(X,\IQ)\right).$$ We denote the
cohomology endowed with this mixed Hodge structure by
$H^*(X,B,\IZ)$. Clearly, the induced rational(!) mixed Hodge
structure is by definition isomorphic to the standard one, but as
soon as $B$ is not integral the decomposition $W_{i+1}=W_i\oplus
\exp(B) H^{i+1}(X,\IQ)$ is not defined over $\IZ$. Also note that
the definition depends on $B$ and not just on $\alpha$, e.g.\ for
$B\in H^{1,1}(X,\IQ)$, which induces the trivial Brauer class,
the mixed Hodge structure will nevertheless be non-split in
general.

ii) In the sequel we shall denote by $H^{*,*}(X,B,\IZ)$
(respectively $H^{*,*}(X,B,\IQ)$) the integral (resp.\ rational)
part of $\exp(B)\left(\bigoplus H^{*,*}(X)\right)$. More precisely,
$H^{p,p}(X,B,\IZ)=\exp(B)(H^{p,p}(X,\IQ))\cap H^*(X,\IZ)$.
Proposition \ref{ProptwistedChern}, iv) then says that
$$\xymatrix{{\rm Im}\left({\rm
ch}^B:K(X,\alpha)\right.\ar[r]&{\phantom{c^B}}
\left.{\phantom{c^B}}\!\!\!\!\!\!\!\!\!\!\!\!\!H^*(X,\IQ)\right)\subset
H^{*,*}(X,B,\IQ)}.$$

iii) Note that the natural product defines homomorphisms of mixed
Hodge structures
$$\xymatrix{H^*(X,B_1,\IZ)\otimes
H^*(X,B_2,\IZ)\ar[r]&H^*(X,B_1+B_2,\IZ)}.$$ All mixed Hodge
structures together yield a mixed Hodge structure on
$H^*(X,\IZ)\otimes_\IZ H^2(X,\IQ)$ that is endowed with an inner
product.
\end{remark}

It seems plausible to generalize the standard Hodge conjecture to
the question whether the map ${\rm ch}^B:K(X,\alpha)\otimes\IQ\to
H^{*,*}(X,B,\IQ)$ is surjective. As it turns out, the twisted
Hodge conjecture is actually equivalent to the untwisted one.
This is the following

\begin{proposition}\label{HC}
Let $X$ be a smooth complex projective variety and $B\in
H^2(X,\IQ)$ a rational B-field inducing the Brauer class
$\alpha$. Then the following two assertions are equivalent:

{\rm i)} The map ${\rm ch}:K(X)_\IQ\to \bigoplus H^{p,p}(X,\IQ)$
is surjective.

{\rm ii)} The map ${\rm ch}^B:K(X,\alpha)_\IQ\to \bigoplus
H^{p,p}(X,B,\IQ)$ is surjective.
\end{proposition}

\begin{proof}
Let us denote the image of ${\rm ch}$ and ${\rm ch}^B$ by $A$
respectively $A_B$. If $E\in K(X,\alpha)$ and $F\in K(X)$, then
$E\otimes F\in K(X,\alpha)$ and ${\rm ch}^B(E\otimes F)={\rm
ch}^B(E)\cdot{\rm ch}(F)$. Hence, $A_B$ is invariant under
multiplication with $A$.

Assume that $A=\bigoplus H^{p,p}(X,\IQ)$, then $A_B\subset
\exp(B)(\bigoplus H^{p,p}(X,\IQ))$ is invariant under
multiplication with $\bigoplus H^{p,p}(X,\IQ)$. As ${\bf
Coh}(X,\alpha)$ contains a locally free sheaf of finite rank,
there exists an element $c\in A_B$ of the form $c=1+{\rm
higher~order~terms}$. From these two statements one easily
deduces that ${\rm ch}^B:K(X,\alpha)_\IQ\to\bigoplus
H^{p,p}(X,B,\IQ)$ is surjective, i.e.\ $A_B=\bigoplus
H^{p,p}(X,B,\IQ)$.

Let us conversely assume that $A_B=\bigoplus H^{p,p}(X,B,\IQ)$.
The idea of the proof is to use the commutativity of the
following diagram
$$\xymatrix{K(X,\alpha)_\IQ\ar[dd]_{(~~)^k}\ar@{>>}[rr]^-{{\rm ch}^B}
&&\bigoplus H^{p,p}(X,B,\IQ)\ar[dd]^{(~~)^k}\\
\\
K(X)_\IQ \ar[rr]^-{{\rm ch}^{B_0}}&&\bigoplus
H^{p,p}(X,B_0,\IQ)}$$ where $B_0:=k\cdot B\in H^2(X,\IZ)$ for
some $k\in\IZ$. As the image of the vertical map on the right
hand side spans, this suffices to conclude that ${\rm ch}^{B_0} $
is surjective as well.

More explicitly, consider  $\beta\in H^{p,p}(X,\IQ)$. Then
$\exp(B)\cdot(1+\beta)\in A_B$ and, therefore, (up to a scaling
factor) $\exp(B)\cdot(1+\beta)={\rm ch}^B(E)-{\rm ch}^B(F)$ for
certain $\alpha$-twisted vector bundles $E$ and $F$. Passing to
the $k$-th power yields $\exp(B_0)\cdot
(1+\beta)^k=\sum(-1)^i{k\choose i}{\rm ch}^{B_0}(E^{\otimes
i}\otimes F^{\otimes {k-i}})=\exp(B_0)\cdot\sum(-1)^i{k\choose
i}{\rm ch}(E^{\otimes i}\otimes F^{\otimes {k-i}})$, because
$E^{\otimes i}\otimes F^{\otimes {k-i}}$ is an untwisted vector
bundle. Thus, $(1+\beta)^k\in A$ and hence $\beta\in A$.
\end{proof}

\begin{remark}
i) There is a simple way to construct a Chern character for
twisted sheaves which only depends on the B-field as an element
in $H^2(X,\IQ)/H^2(X,\IZ)$. Indeed, one may consider
$\exp(-B){\rm ch}^B(~~)$. Changing $B$ by an integral B-field
does not affect this expression. However, the above approach has
the advantage that the denominators that occur in ${\rm ch}^B\in
H^*(X,\IQ)$ are universal. This fact will be important for the
construction of certain Hodge isometries over $\IZ$ (and not only over
$\IQ$), in Section \ref{TwistedFMonCohSection}.

ii) Twisted Chern characters are alluded to at different places
in the lite\-ra\-ture, but we couldn't find the above explicit
construction. However, Eyal Markman informed us that Jun Li has
developed a theory of of connections on twisted holomorphic
bundles and a twisted analogue of the Hermite--Einstein equation.
\end{remark}
\section{Generalized CY structures versus twisted K3
surfaces}\label{GCYSection}

In this section we shall compare the notions of  twisted K3
surfaces and generalized  Calabi--Yau structures on $M$.

\begin{definition}\label{DeftwistedK3}
A \emph{twisted K3 surface} $(X,\alpha)$ consists of a K3 surface
$X$ together with a Brauer class $\alpha\in{\rm Br}(X)$. We say
that $(X,\alpha)\cong (X',\alpha')$ if there exists an isomorphism
$f:X\cong X'$ with $f^*\alpha'=\alpha$.
\end{definition}

\begin{definition}
A \emph{generalized Calabi--Yau} structure on $M$ is an even
closed complex form $\varphi=\varphi_0+\varphi_2+\varphi_4$ such
that $\langle\varphi,\varphi\rangle=0$ and
$\langle\varphi,\overline\varphi\rangle>0$.
\end{definition}

Here, $\langle~~,~~\rangle $ is the Mukai pairing on the level of
forms and the positivity of
$\langle\varphi,\overline\varphi\rangle$ is meant with respect to
a fixed volume form. For details see \cite{HuyGen}. In this paper
we are only interested in generalized Calabi--Yau structures of
the form $\varphi=\exp(B)\cdot \sigma= \sigma+B\wedge\sigma$,
where $\sigma$ is a holomorphic two-form on a K3 surface
$X=(M,I)$ and $B$ a real closed two-form. If $B$ can be chosen
such that $[B]$ is a rational class, then $\varphi$ is called a
\emph{rational generalized Calabi--Yau structure}. In the sequel,
we often just work with the cohomology classes of $\varphi$ and
$B$, which for simplicity will be denoted by the same symbols

Suppose $X=(M,I)$ is a given K3 surface. To any rational B-field
$B\in H^2(M,\IQ)$ one can associate the twisted K3 surface
$(X,\alpha_B)$ (for the defi\-ni\-tion of $\alpha_B$ see the
introduction) and a generalized Calabi--Yau structure
$\varphi=\exp(B)\cdot \sigma$.

In fact, $(X,\alpha_B)$ depends only on the cohomology class of
$\varphi$ and we therefore get a natural map
$$\xymatrix{\Big\{{\rm rational~gen.~CYs~}
[\exp(B)\cdot \sigma]\Big\}\ar[r]&\Big\{{\rm
twisted~K3s~}(X,\alpha)\Big\},}$$
 which is surjective due to $H^3(X,\IZ)=0$.
\begin{definition}
Let $X$ be a K3 surface with a (rational) B-field $B\in
H^2(X,\IQ)$. Then we denote by $\widetilde H(X,B,\IZ)$ the
weight-two Hodge structure on $H^*(X,\IZ)$ with
$$\widetilde H^{2,0}(X,B):=\exp(B)\left(H^{2,0}(X)\right)$$
and $\widetilde H^{1,1}(X,B)$ its orthogonal (with respect to the
Mukai pairing) complement.
\end{definition}

Clearly, $\widetilde H(X,B,\IZ)$ only depends on the generalized
Calabi--Yau structure $\varphi=\sigma+B\wedge\sigma$. In other
words, $B$ and $B'$ define the same weight-two Hodge structure on
$H^*(M,\IZ)$ if and only if $B^{0,2}={B'}^{0,2}\in H^2(X,\ko_X)$.
If $B$ and $B'$ differ by an integral class $B_0\in H^2(M,\IZ)$
then
$$\xymatrix{\exp(B_0):\widetilde H(X,B,\IZ)\ar[r]&\widetilde
H(X,B',\IZ)}$$ is a Hodge isometry. This yields the diagram
$$\xymatrix{\Big\{{\rm rational~gen.~CYs~}
[\exp(B)\cdot \sigma]\Big\}\ar[d]\ar[r]&\Big\{{\rm
twisted~K3s~}(X,\alpha)\Big\}\ar[d]\\
\Big\{{\rm weight}-2~{\rm HS~on~}H^*(M,\IZ)\Big\}\ar[r]&\Big\{{\rm
weight}-2 ~{\rm HS~on~}H^*(M,\IZ)\Big\}/_{\cong}}$$

This new Hodge structure comes along with a natural orientation
of its positive directions.
We shall briefly explain what this means.
If $X$ is a K3 surface with $\sigma$ a
generator of $H^{2,0}(X)$ and $\omega$ a K{\"a}hler class (e.g.\ an
ample class if $X$ is algebraic), then $\langle{\rm Re}(\sigma),
{\rm Im}(\sigma),1-\omega^2/2,\omega\rangle$ is a positive
four-space in $\widetilde H(X,\IR)$ which comes, by the choice of
the basis, with a natural orientation. It is easy to see that this
orientation is independent of the choice of $\sigma$ and $\omega$.
Let $g:\Gamma\to \Gamma'$ be an isometry of lattices with
signature $(4,t)$. Suppose positive four-spaces
$V\subset\Gamma_\IR$ and $V'\subset\Gamma'_\IR$ and orientations
for both of them have been chosen. Then one says that $g$
preserves the given orientation of the positive directions
(or, simply, that $g$ is orientation preserving) if the composition of
$g_{_\IR}:V\to\Gamma'_\IR$ and the orthogonal projection
$\Gamma'_\IR\to V'$ is compatible with the given orientations of
$V$ and $V'$. By $\OO_+(\Gamma)$ one denotes the group of all
orientation preserving  orthogonal transformations.
We shall also use the analogous notation $\OO_+(U)$ for the hyperbolic
plane $U=H^0\oplus H^4$.

In \cite{HuyGen} we explained how to associate to any generalized
Calabi--Yau structure $\varphi\in H^*(M,\IC)$ the
\emph{generalized (or twisted) Picard group} and the
\emph{generalized~(or twisted) transcendental lattice}. Let us
discuss these two lattices a bit further in the case of a
generalized Calabi--Yau structure $\varphi=\exp(B)\cdot\sigma$
with $\sigma$ the holomorphic two-form on the K3 surface
$X=(M,I)$. By definition the generalized Picard group is
$${\rm Pic}(X,\varphi):=\{\beta\in
H^*(M,\IZ)~|~\langle\beta,\varphi\rangle=0\}$$ and the
\emph{generalized~transcendental lattice}
$$T(X,\varphi):={\rm Pic}(X,\varphi)^\perp\subset H^*(M,\IZ),$$
where the orthogonal complement is taken with respect to the
Mukai pairing. The latter comes with a natural weight-two Hodge
structure which will always be understood. Note that both
lattices depend only on the Hodge structure $\widetilde
H(X,B,\IZ)$ and, therefore, their isomorphism classes only on the
induced Brauer class $\alpha_B$. Later we often write ${\rm
Pic}(X,B)$ and $T(X,B)$ or, if only the isomorphism type is
relevant, ${\rm Pic}(X,\alpha_B)$ and $T(X,\alpha_B)$.

By definition, ${\rm Pic}(X,\varphi)$ is the set of integral
classes $(\delta_0,\delta_2,\delta_4)$ with
$\delta_0\int\sigma\wedge B=\int\sigma\wedge\delta_2$. In
particular, no condition on $\delta_4$ and, hence, ${\rm
Pic}(X)\oplus H^4(X,\IZ)\subset{\rm Pic}(X,\varphi)$. One also
has the following description of the generalized Picard group
$${\rm Pic}(X,B)={\rm Pic}(X,\varphi)=H^{*,*}(X,B,\IZ)=\widetilde
H^{1,1}(X,B,\IZ).$$

Indeed, if one writes
$(\delta_0,\delta_2,\delta_4)=\exp(B)(\alpha_0,\alpha_2,\alpha_4)$,
one finds that $\alpha_2\in H^{1,1}(X)$, i.e.\
$\int\sigma\wedge\alpha=0$, is equivalent to
$\delta_0\int\sigma\wedge B=\int\sigma\wedge\delta_2$. Hence,
${\rm Pic}(X,\varphi)=H^{*,*}(X,B,\IZ)$.

The third equality follows from the fact that $\exp(B)$ is an
orthogonal transformation with respect to the Mukai pairing,
i.e.\ $\langle\exp(B)(~~),\exp(B)(~~)\rangle=\langle~~,~~\rangle$.

The next result can be seen as a consequence of Proposition
\ref{HC}.

\begin{corollary}
Let $X$ be a K3 surface and $\alpha:=\exp(B^{0,2})\in
H^2(X,\ko_X^*)$ a given Brauer class induced by a rational B-field
$B$. Then the twisted Chern character
$$\xymatrix{{\rm ch}^B:K(X,\alpha)\ar@{>>}[r]&{\rm Pic}(X,B)}$$
is surjective.\qqed
\end{corollary}


It is easy to see that there exists a finite index immersion
$$\xymatrix{{\rm Pic}(X)\oplus (\lambda_Bu_2+\lambda_BB)\IZ\oplus
H^4(X,\IZ)\ar@{^({-}>}[r]&{\rm Pic}(X,B)}$$ for a certain integer
$\lambda_B$ such that $\lambda_BB\in H^2(X,\IZ)$. In particular, ${\rm
Pic}(X,B)$ and ${\rm Pic}(X)\oplus U$ are always of the same
rank. Moreover, $X$ is algebraic if and only if ${\rm Pic}(X,B)$
contains two positive directions.

We also need to clarify the relation  between $T(X,B)$ and the
transcendental lattice $T(X)$. In \cite{HuyGen} we have shown
that $\exp(-B)$ defines a Hodge isometry $T(X,B)\cong
T(X,\alpha_B)$, where $T(X,\alpha_B)$ is the kernel of the natural
map $(B,~~):T(X)\to\IQ/\IZ$ defined by the intersection product
with $B$ (it only depends on $\alpha_B$). Once more, $T(X,B)$ and
$T(X)$ are therefore lattices of the same rank.


\section{Some remarks on Brauer groups of K3
surfaces}\label{BrauerSect}

By definition, the Brauer group ${\rm Br}(X)$ of $X$ is the set
of torsion classes in $H^2(X,\ko_X^*)$. If the order of a class
$\alpha\in{\rm Br}(X)$ divides $k$ then there exists a B-field
lift $B$ of $\alpha$, i.e.\ a class $B\in H^2(X,\IQ)$ with
$\alpha_B=\alpha$, such that $kB\in H^2(X,\IZ)$. Using the
intersection product with $B$ yields a linear map
$(B,~~):T(X)\to\IQ$ with image contained in $\frac{1}{k}\IZ$ and,
by dividing by $\IZ\subset \frac{1}{k}\IZ$, a linear map
$\alpha:T(X)\to \frac{1}{k}\IZ/\IZ\cong\IZ/k\IZ$. It is
straightforward to check that this construction yields an
isomorphism $${\rm Br}_{k-{\rm tor}}(X)\cong\Hom(T(X),\IZ/k\IZ).$$
Moreover, if the order of $\alpha$ equals $k$ then the induced
map $\alpha:T(X)\to\IZ/k\IZ$ is surjective.

If one is not interested in any specific torsion, the
construction yields an alternative description of the Brauer
group as ${\rm Br}(X)=\Hom(T(X),\IQ/\IZ)$.

\begin{remark}\label{ordertrans}
 In the last paragraph, we used the notation $T(X,\alpha)$ for
the kernel of the map $T(X)\to\IQ/\IZ$. Thus, the order of a
Brauer class $\alpha$ can be computed by using the standard
formula
$$|\alpha|^2\cdot\left|{\rm disc}( T(X,\alpha))\right|=\left|{\rm
disc}(T(X))\right|.$$

Moreover, the existence of the aforementioned Hodge isometry
$\exp(-B):T(X,B)\cong T(X,\alpha)$ proves that the order of a
Brauer class is encoded by the twisted transcendental lattice
$T(X,B)$, where $B$ is an arbitrary rational B-field lift of
$\alpha$ and, of course, the standard transcendental lattice
$T(X)$.
\end{remark}

But even with this description the Brauer group ${\rm Br}(X)$,
which is abstractly isomorphic to $(\IQ/\IZ)^{22-\rho(X)}$, is
otherwise a rather mysterious object. In order to get a better
understanding of it, we shall introduce various equivalence
relations.

\begin{definition}\label{EquivDef}
We define \emph{derived equivalence} of two Brauer classes
$\alpha,\alpha'\in{\rm Br}(X)$ by
$$\alpha\stackrel{D}{\sim}\alpha' \Longleftrightarrow \Db(X,\alpha)\cong
\Db(X,\alpha')$$ where the isomorphism on the right hand side
means Fourier--Mukai equivalence.

\emph{Hodge equivalence} is defined by
$$\alpha\stackrel{H}{\sim}\alpha' \Longleftrightarrow \widetilde H(X,B,\IZ)\cong
\widetilde H(X,B',\IZ),$$ where the isomorphism on the right hand
side means Hodge isometry.

We say that two Brauer classes $\alpha,\alpha'$ are
\emph{$T$-equivalent}, $\alpha\stackrel{T}{\sim}\alpha'$ if and
only if there exists a Hodge isometry $T(X,\alpha)\cong
T(X,\alpha')$.
\end{definition}

The Hodge structures $\widetilde H(X,B,\IZ)$ and $\widetilde
H(X,B',\IZ)$ are defined in terms of $B,B'\in H^2(X,\IQ)$, which
are chosen such that $\alpha=\alpha_B$ and $\alpha'=\alpha_{B'}$.
Note that Hodge equivalence is well-defined, as different choices
of the B-field lift of a Brauer class give rise to isomorphic
Hodge structures. Derived equivalences will only be considered
for algebraic K3 surfaces and we will assume that the equivalence
is a Fourier--Mukai equivalence. Hodge equivalence and
$T$-equivalence of Brauer classes make perfect sense also in the
non-algebraic situation.

\begin{remark}\label{implications}
i) Clearly, two Hodge equivalent Brauer classes are
also $T$-equivalent. Building upon the construction of the last
section, we shall show that  derived equivalence implies  Hodge
equivalence (see Corollary \ref{H=D}). So eventually, we will have
$$\alpha\stackrel{D}{\sim}\alpha'\Longrightarrow
\alpha\stackrel{H}{\sim}\alpha'\Longrightarrow\alpha\stackrel{T}{\sim}\alpha'.$$
The first implication is expected to be almost an equivalence (see
Conjecture \ref{Calconjref}). More precisely, if
$\alpha\stackrel{H}{\sim}\alpha'$ such that the Hodge isometry in
Definition
\ref{EquivDef} can be chosen orientation preserving then
  $\alpha\stackrel{D}{\sim}\alpha'$.

ii) Note that by Remark \ref{ordertrans} $T$-equivalent Brauer
classes are of the same order. This immediately shows that there
always exists an infinite number of $T$-equivalence classes. The
same holds true for Hodge equivalence (see i)) and for derived
equivalence (see Corollary \ref{infintcor}).
\end{remark}

The automorphism group ${\rm Aut}(X)$ acts naturally on ${\rm
Br}(X)$ via $\alpha\mapsto f^*\alpha$. Note that
$f^*\alpha\stackrel{*}{\sim}\alpha$ for $*=D,H,T$. Indeed, the
induced action of $f$ on cohomology yields a Hodge isometry
between the two Hodge structures, which proves the equivalence
for $*=H,T$, and $E\mapsto f^*E$ defines an equivalence
$\Db(X,\alpha)\cong\Db(X,f^*\alpha)$. In other words, ${\rm
Aut}(X)$ preserves the equivalence classes of all three
equivalence relations.

 Since the automorphism group
of a K3 surface might very well be infinite, the set
$$\{\alpha\in{\rm Br}(X)~|~\alpha\stackrel{*}{\sim}\alpha_0\}$$  with  $\alpha_0\in{\rm Br}(X)$
fixed will in general be infinite. Note however that for
algebraic K3 surfaces the action of ${\rm Aut}(X)$ on
$H^2(X,\ko_X)$ and hence on $H^2(X,\ko_X^*)$ factorizes over a
finite group (see the proof of the following proposition).

\begin{proposition}\label{propfinitebrauer}
Let $X$ be a K3 surface (not necessarily algebraic) and
$\alpha_0\in{\rm Br}(X)$. Then
$$\raisebox{.0ex}{$\{\alpha~|~\alpha\stackrel{*}{\sim}\alpha_0\}$}\raisebox{-1.0ex}{$/$}\raisebox{-1.30ex}{${\rm Aut}(X)$}$$
is finite. Moreover, if $X$ is algebraic then
$\{\alpha~|~\alpha\stackrel{*}{\sim}\alpha_0\}$ is finite. In
both cases, $*=H$ or $*=T$.
\end{proposition}

\begin{proof}
Since Hodge equivalence implies $T$-equivalence, it suffices to
prove the assertions for $*=T$. The second follows easily from
the first, as the action of the possibly infinite group ${\rm
Aut}(X)$ on $H^2(X,\ko_X)$ factorizes over a finite group if $X$
is algebraic (see the comment at the end of the proof).

Let us suppose $\alpha_0$ is of order $k$. Then the same holds
true for any $\alpha\stackrel{T}{\sim}\alpha_0$ (see Remark
\ref{implications}, ii)). For any such $\alpha$ we fix a Hodge
isometry $T(X,\alpha)\cong T(X,\alpha_0)=:T_0$.

As was explained earlier, for any Brauer class $\alpha$ of order
$k$ there exists a natural short exact sequence
$$\xymatrix{0\ar[r]&T(X,\alpha)\ar[r]&T(X)\ar[r]&\IZ/k\IZ\ar[r]&0.}$$
Moreover, the class $\alpha$ itself is determined by the map
$T(X)\to\IZ/k\IZ$.

Therefore, any class $\alpha\stackrel{T}{\sim}\alpha_0$ is
determined by a Hodge embedding $T_0\hookrightarrow T(X)$
compatible with the intersection form  and an isomorphism
$T(X)/T_0\cong\IZ/k\IZ$.

Clearly, for a given embedding $T_0\hookrightarrow T(X)$ there is
only a finite number of isomorphisms $T(X)/T_0\cong\IZ/k\IZ$.
Thus, it suffices to show that the number of embeddings
$T_0\hookrightarrow T(X)$ that are compatible with Hodge
structure and intersection form is finite up to the action of
${\rm Aut}(X)$.

Up to the action of the group of isometries $\OO(T(X))$ the
embedding $T_0\hookrightarrow T(X)$ is determined by a subgroup
of the finite group $T_0\dual/T_0$. Hence, up to the action of
$\OO(T(X))$, there is only a finite number of possibilities for
$T_0\hookrightarrow T(X)$. On the other hand, any isometry
$g\in\OO(T(X))$ fixing $T_0$ is in fact a Hodge isometry of
$T(X)$.

To conclude, one uses the fact that the image of the natural map
${\rm Aut}(X)\to {\rm Aut}(T(X))$ is a finite index subgroup.
Here, ${\rm Aut}(T(X))$ is the group of Hodge isometries of the
transcendental lattice $T(X)$.

For the reader's convenience we include a proof of this fact. Up
to the action of $\OO(\Lambda)$ there exists only finitely many
primitive embeddings $T(X)\hookrightarrow\Lambda$. In particular,
up to isometries of $\Lambda$ there exist only  finitely many
Hodge isometries $T(X)\cong T(X)$. Now use that any
$g\in\OO(\Lambda)$ that is compatible with a Hodge isometry of
$T(X)$, is in fact a Hodge isometry $H^2(X,\IZ)$. By the global
Torelli theorem one knows that up to sign any Hodge isometry of
$H^2(X,\IZ)$ modulo the action of the Weyl group, which acts
trivially on the transcendental lattice $T(X)$, is induced by an
automorphism of $X$. This proves the assertion.

Note that in the case of an algebraic K3 surface, i.e.\ when the
signature of $T(X)$ is $(2,s)$,  the group ${\rm Aut}(T(X))$ is
isomorphic to a discrete subgroup of the compact group
$\OO(2)\times\OO(s)$. Hence, the action of ${\rm Aut}(X)$ on
$H^2(X,\ko_X)$ factorizes over a finite group.\end{proof}


In Section \ref{TwistedFMonCohSection} we shall see that the
proposition implies the analogous statement for derived
equivalence.


\section{Twisted Fourier--Mukai equivalence on
cohomology}\label{TwistedFMonCohSection}

We indicate how to modify the well-known arguments in order to
make them work in the twisted situation. We start with an
arbitrary smooth projective variety and shall restrict to
algebraic K3 surfaces later on.

\begin{definition}
The \emph{Mukai vector} of $E\in K(X,\alpha)$ is
$$v^B(E):={\rm ch}^B(E)\cdot\sqrt{{\rm td}(X)}.$$
\end{definition}

As before, $B\in H^2(X,\IQ)$ is a rational B-field such that
$\exp(B^{0,2})=\alpha$.

Clearly, $$\xymatrix{v^B:K(X,\alpha)\ar[r]& H^*(X,\IQ)}$$ takes
again values in $\exp(B)\left(\bigoplus H^{p,p}(X)\right)$ (see
Remark \ref{MHSRem}, ii)).

\begin{remark}
With this definition the Riemann--Roch formula still holds, i.e.\
for $E,F\in{\bf Coh}(X,\alpha)$ one has
$$\chi(E,F)=\langle v^B(E)\dual. v^B(F)\rangle,$$
where $\langle~~,~~\rangle$ is the Mukai pairing (or rather its
generalization introduced by C\u{a}ld\u{a}raru \cite{Cal3}). This
follows easily from the observations
$\chi(E,F)=\chi(X,E\dual\otimes F)$, ${\rm ch}^B(E)\dual={\rm
ch}^{-B}(E\dual)$, and ${\rm ch}^{-B}(E\dual).{\rm ch}^B(F)={\rm
ch}(E\dual\otimes F)$.
\end{remark}

Let $X$ and $X'$ be two smooth projective varieties equipped with
topolo\-gi\-cally trivial Brauer classes $\alpha$ respectively
$\alpha'$ and B-field lifts $B$ respectively $B'$.

We then consider the natural B-field $(-B)\boxplus
B':=q^*(-B)+p^*B'\in H^2(X\times X',\IQ)$ and the induced Brauer
class $\alpha^{-1}\boxtimes\alpha'\in H^2(X\times X',\ko^*)$.

Any $e\in K(X\times X',\alpha^{-1}\boxtimes\alpha')$ defines a
Fourier--Mukai transformation

$$\xymatrix{\Phi^K_e:K(X,\alpha)\ar[r]&K(X',\alpha').}$$

As in the untwisted case, we obtain a commutative diagram

$$\xymatrix{K(X,\alpha)\ar[d]_{v^B}\ar[rr]^{\Phi^K_e}&&K(X',\alpha')\ar[d]^{v^{B'}}\\
H^*(X,\IQ)\ar[rr]_-{\Phi^H_{v^{(-B)\boxplus
B'}(e)}}&&H^*(X',\IQ)}$$

The idea here is to write ${\rm ch}^B(E)$ as ${\rm ch}(E_B)$ as in
the proof of Proposition \ref{ProptwistedChern}.
The original argument uses the
Grothendieck--Riemann--Roch formula at this point.
Fortunately, due to work of Atiyah and Hirzebruch,
see \cite[Thm.\ 1, Thm.\ 2]{AH}, the same formula holds true also in the
differentiable setting.

\bigskip

Suppose a Fourier--Mukai transform
$$\xymatrix{\Phi:\Db(X,\alpha)\ar[r]&\Db(X',\alpha')}$$
 with kernel $\ke\in\Db(X\times
X',\alpha^{-1}\boxtimes\alpha')$ is given (see \cite{Cal} for the
discussion of all the necessary derived functors). We will always
assume that a rational B-field has been chosen inducing a given
Brauer class $\alpha=\alpha_B$.

Many of the following arguments are taken directly from the
original sources \cite{Cal,Mu}. We only indicate the necessary
modifications.

\medskip

1) {\it The right and left adjoints of a Fourier--Mukai functor
$\Phi:\Db(X,\alpha)\to\Db(X',\alpha')$ exist. They are of
Fourier--Mukai type.}

\medskip

2) {\it The composition of two twisted Fourier--Mukai functors is
again of Fourier--Mukai type.}

\medskip

3) {\it If a Fourier--Mukai equivalence
$\Phi:\Db(X,\alpha)\to\Db(X,\alpha)$ with kernel $\ke$ is
isomorphic to the identity, then $\ke\cong\ko_{\Delta}$.}
\mynote{Check} Observe that $\ko_{\Delta}$ is indeed an object in
$\Db(X\times X,\alpha^{-1}\boxtimes\alpha)$.

\medskip

4) {\it To any Fourier--Mukai functor
$\Phi:\Db(X,\alpha)\to\Db(X',\alpha')$ with kernel $\ke$ one
associates a homomorphism $\Phi_*^{B,B'}:H^*(X,\IQ)\to
H^*(X',\IQ)$ depending on the B-field lifts of $\alpha$ and
$\alpha'$. If $\Phi$ is a Fourier--Mukai equivalence, then
$\Phi_*^{B,B'}$ is bijective.}

Here, the cohomological Fourier--Mukai transform is with respect
to the Mukai vector $v^{(-B)\boxplus B'}(\ke)$. Check that
$v^{(-B)\boxplus B}(\ko_\Delta)=[\Delta]$. In fact, a little more
can be said. As in the untwisted case, a Fourier--Mukai
equivalence  induces a cohomological Fourier--Mukai transform
which yields  isomorphisms $\bigoplus_{p-q=a}
H^{p,q}(X,B,\IQ)\cong\bigoplus_{p-q=a} H^{p,q}(X',B',\IQ)$ for all
$a$.

\medskip

5) {\it Any Fourier--Mukai equivalence
$\Phi:\Db(X,\alpha)\to\Db(X,\alpha)$ induces a rational isometry
$\Phi_*^{B,B'}$ with respect to C\u{a}ld\u{a}raru's
generalization of the Mukai pairing.}

\medskip

So far, $X$ was just any smooth projective variety. We now
restrict to the case of K3 surfaces. We shall use the notation of
the last section.

\medskip

 6) {\it Let $X$ and $X'$ be two K3 surfaces and $e\in
K(X\times X',\alpha^{-1}\boxtimes\alpha')$. Then $v^{(-B)\boxplus
B'}(e)\in H^*(X\times X',\IZ)$.}

Once more, one has to replace the Grothendieck--Riemann--Roch
formula in Mukai's original proof by the differentiable version
of it (cf.\ \cite{AH}).

Here comes the main result of this section. It is an immediate
consequence of the above results.

\begin{proposition}\label{HyieldsD}
Let $X$ and $X'$ be two algebraic K3 surfaces equipped with
rational B-fields $B$ respectively $B'$. The induced elements in
the Brauer group are denoted $\alpha$ respectively $\alpha'$. If
$$\Phi:\Db(X,\alpha)\cong\Db(X',\alpha')$$
is a Fourier--Mukai equivalence, then the induced cohomological
Fourier--Mukai map
$$\Phi_*^{B,B'}:\widetilde H(X,B,\IZ)\cong \widetilde H(X',B',\IZ)$$
is an isometry of integral weight-two Hodge structures.\qqed
\end{proposition}

\begin{corollary}\label{H=D}
Let $\alpha,\alpha'\in{\rm Br}(X)$. Then
$\alpha\stackrel{D}{\sim}\alpha'$ implies
$\alpha\stackrel{H}{\sim}\alpha'$.\qqed
\end{corollary}

Together with Proposition \ref{propfinitebrauer} this also yields:
\begin{corollary}\label{finiteBrmodD}
Let $X$ be an algebraic K3 surface  and $\alpha_0\in{\rm Br}(X)$.
Then the equivalence class
$\{\alpha~|~\alpha\stackrel{D}{\sim}\alpha_0\}$ is a finite
set.\qqed
\end{corollary}

The corollary is in fact a special case of the more general
statement that up to automorphism the number of Fourier--Mukai
partners is finite. Here a twisted K3 surface $(X',\alpha')$ is a
\emph{(twisted) Fourier--Mukai partner} of a given twisted K3
surface $(X,\alpha)$ if there exists a Fourier--Mukai equivalence
$\Db(X,\alpha)\cong\Db(X',\alpha')$. In \cite{BM} it has been
shown that for any given untwisted K3 surface there are only
finitely many    untwisted Fourier--Mukai partners. As a
consequence of Proposition \ref{HyieldsD} and its Corollary
\ref{finiteBrmodD} one deduces in the same way the finiteness in
the twisted case:

\begin{corollary}\label{finiteFMpartners}
Any twisted algebraic K3 surface $(X,\alpha)$ admits only
finitely many Fourier--Mukai partners up to isomorphisms.
\end{corollary}

\begin{proof}
Any Fourier--Mukai equivalence
$\Db(X,\alpha_B)\cong\Db(X',\alpha_{B'})$ induces a Hodge
isometry $T(X,B)\cong T(X',B')$. As was explained earlier, the
transcendental lattice $T(X,B)$ can be embedded  into $T(X)$ via
$\exp(-B)$ and similarly for $X'$. Hence $T(X')$ sits between
$T(X,B)$ and $T(X,B)\dual$, i.e.\
$$T(X,B)\cong T(X',B')\subset T(X')\subset T(X',B')\dual\cong
T(X,B)\dual.$$ This shows that there are only finitely many
 possibilities for the isomorphism type of $T(X')$. Arguing as in
\cite{BM} shows that there exist only  finitely many isomorphism
classes of K3 surfaces $X'$ for which a Brauer class $\alpha'$
can be chosen such that $\Db(X,\alpha_B)$ and $\Db(X',\alpha')$
are Fourier--Mukai equivalent.

Eventually, Corollary \ref{finiteBrmodD} says that on any of the
finitely many K3 surfaces $X'$ one only has a finite number (up to
isomorphisms) of Brauer classes realizing a derived category that
is Fourier--Mukai equivalent to the given one $\Db(X,\alpha)$.
\end{proof}

As another consequence of the proposition and Remark
\ref{implications}, ii) one obtains

\begin{corollary}\label{infintcor}
Let $X$ be an algebraic K3 surface. Then there exists an infinite
number of pairwise inequivalent twisted derived categories
$\Db(X,\alpha)$.\qqed
\end{corollary}

We conclude this section with a detailed discussion of
C\u{a}ld\u{a}raru's conjecture. It turns out that the original
formulation has to be modified. At the same time, we shall
propose a version that relates Hodge isometries and derived
equivalences in a more precise way.

Let us begin with the untwisted version. We shall state it as a
conjecture, although most, but not all, of it has been proved
already.

\begin{conjecture}
Let $X$ and $X'$ be two algebraic K3 surfaces.

{\rm a)} If $\Phi:\Db(X)\cong\Db(X')$ is a Fourier--Mukai
equivalence, then $\Phi_*:\widetilde H(X,\IZ)\to \widetilde H
(X',\IZ)$ satisfies:

\begin{itemize}
\item[i)] $\Phi_*$ is a Hodge isometry and
\item[ii)] $\Phi_*$ preserves the natural orientation of the four positive
directions.
\end{itemize}

{\rm b)} If $g:\widetilde H(X,\IZ)\to \widetilde H (X',\IZ)$
satisfies {\rm i)} and {\rm ii)} then there exists a
Fourier--Mukai equivalence $\Phi$ with $g=\Phi_*$.
\end{conjecture}

The known results are essentially due to Mukai and Orlov with
complementary observations provided by \cite{HLOY1,Ploog}. That
the orientation should be important at this point was first
observed in \cite{Sz}. Part b) of the conjecture is essentially
known.  The only missing detail that we could not find in the
literature is explained in the next section.
In a) one only knows, for the time being, that
$\Phi_*$ satisfies i). The problem of showing ii) can be reduced
to the question whether the Hodge isometry $j=(-{\rm
id}_{H^0})\oplus {\rm id}_{H^2}\oplus(-{\rm id}_{H^4})$ can be
realized by an autoequivalence.

Let us now state the twisted version.

\begin{conjecture}\label{Calconjref}
Let $X$ and $X'$ be two algebraic K3 surfaces with B-fields $B$
and $B'$.

{\rm a)} If $\Phi:\Db(X,\alpha_B)\cong\Db(X',\alpha_{B'})$ is a
Fourier--Mukai equivalence, then $\Phi_*^{B,B'}:\widetilde
H(X,B,\IZ)\to \widetilde H (X',B',\IZ)$ satisfies:

\begin{itemize}
\item[i)] $\Phi_*^{B,B'}$ is a Hodge isometry and
\item[ii)] $\Phi_*^{B,B'}$ preserves the natural orientation of the four positive
directions.
\end{itemize}

{\rm b)} If $g:\widetilde H(X,B,\IZ)\to \widetilde H (X',B',\IZ)$
satisfies {\rm i)} and {\rm ii)} then there exists a
Fourier--Mukai equivalence $\Phi$ with $g=\Phi_*^{B,B'}$.
\end{conjecture}

Proposition \ref{HyieldsD} shows that i)  part a) holds. Later we
will explain how to deduce part b) for large Picard number.

\begin{remark}
C\u{a}ld\u{a}raru stated his conjecture originally as:
$\Db(X,\alpha)\cong \Db(X',\alpha')$ if and only if there exists
a Hodge isometry $T(X,\alpha)\cong T(X',\alpha')$.

It seems two problems may occur. First of all, unlike the
untwisted case the existence of a Hodge isometry of the twisted
transcendental lattices $T(X,B)\cong T(X,\alpha_B)\cong
T(X',\alpha_{B'})\cong T(X,B)$ does not necessarily yield the
existence of a Hodge isometry $\widetilde H(X,B,\IZ)\cong
\widetilde H(X',B',\IZ)$. Indeed, the orthogonal complement of
the twisted transcendental lattice $T(X,B)$ does not, in general,
contain a hyperbolic plane, so that Nikulin's results do not apply
anymore. See Example \ref{counter} for a counterexample.

Secondly, even when a Hodge isometry $\widetilde H(X,B,\IZ)\cong
\widetilde H(X',B',\IZ)$ can be found it might reverse the
orientation of the positive directions. In the untwisted case,
this is no problem, as one might compose with $-j$ in order to
get a Hodge isometry that preserves the orientation of the
positive directions. This trick does not work any longer in the
twisted case, as $-j$ is a Hodge isometry between $\widetilde
H(X,B,\IZ)$ and $\widetilde H(X,-B,\IZ)$.

On the level of twisted derived categories this is related to the
problem  that $\Db(X,\alpha)$ and $\Db(X,\alpha^{-1})$ might a
priori not be Fourier--Mukai equivalent if $\alpha$ is of order
$>2$ (compare with Corollary \ref{dualtwisted}).
\end{remark}

\begin{example}\label{counter}
We shall construct a twisted K3 surface $(X,\alpha=\alpha_B)$ with
$\alpha$ of order five such that $T(X,\alpha)\cong T(X,\alpha^2)$,
but the two twisted Hodge structures $\widetilde H(X,B,\IZ)$ and
$\widetilde H(X,2B,\IZ)$ are not Hodge isometric. The latter will
be ensured by showing that ${\rm Pic}(X,B)$ and ${\rm Pic}(X,2B)$
are not isometric.

To be very explicit we write $\Lambda$ as $\Lambda\cong U_1\oplus
\Lambda'$, with $U_1$ isomorphic to the hyperbolic plane. The
standard basis is denoted $e_1,e_2$. Let $X$ be a K3 surface with
maximal transcendental lattice $T(X)=\Lambda$ and choose
$B=(1/5)(e_1+e_2)$. It is not difficult to check that $\alpha_B$
under these assumptions is indeed of order five. The twisted
Picard groups ${\rm Pic}(X,B)$ and ${\rm Pic}(X,2B)$ can now be
described as follows ${\rm Pic}(X,B)=u_1\IZ\oplus
(5u_2+(e_1+e_2))\IZ$ respectively ${\rm Pic}(X,2B)=u_1\IZ\oplus
(5u_2+2(e_1+e_2))\IZ$ (see Section \ref{GCYSection}). As abstract
lattices they are given by the matrices
$$\left(\begin{array}{cc}
0&-5\\
-5&2
\end{array}\right)~~{\rm respectively}~~\left(\begin{array}{cc}
0&-5\\
-5&8
\end{array}\right).$$
These two matrices define inequivalent lattices, for $2$ is
realized as $((0,1),(0,1))_1$ in the first, but there is no
integral vector $(a,b)$ such that $((a,b),(a,b))_2=2$.

Note that the K3 surface $X$ is not algebraic, but it is possible
to start with this example and construct an algebraic
counterexample by adding an algebraic Picard group of the form
$40\IZ$. Indeed, the two lattices given by the symmetric matrices
$$\left(\begin{array}{ccc}
0&-5&0\\
-5&2&0\\
0&0&40
\end{array}\right)~~{\rm respectively}~~\left(\begin{array}{ccc}
0&-5&0\\
-5&8&0\\
0&0&40
\end{array}\right)$$
are not equivalent.

In order to see this, it suffices to verify that
$-10xy+8y^2+40z^2=2$ or, equivalently, $-5xy+4y^2+20z^2=1$ has no
integral solution.

Suppose $(x,y,z)$ is an integral solution. Then $x$ and $y$ are
both odd and hence $p:=(x-y)/2$, $q:=(5x-3y)/2$ are integral.
Replacing $x=q-3p$ and $y=q-5p$ in the above equation yields
$25p^2-1=q^2-20z^2$. Viewing this equation modulo four reveals
that $p$ is odd and $q$ is even. We write $p=2a+1$ and $q=2b$
with $a,b\in\IZ$. Thus, $(5a+2)(5a+3)=b^2-5z^2$.

Let $d:=(b,z)$ and write $b=d\beta$, $z=d\zeta$. Then
$d^2|(5a+2)$ or $d^2|(5a+3)$ and, therefore, $mn=\beta^2-5\zeta^2$
with either $5a+2=md^2$, $5a+3=n$ or, respectively, $5a+2=m$,
$5a+3=nd^2$. As $(\zeta,m)=1=(\zeta,n)$, this shows
$\left(\frac{5}{m}\right)=1=\left(\frac{5}{n}\right)$. On the
other hand, one has
$\left(\frac{m}{5}\right)=-1=\left(\frac{n}{5}\right)$. Indeed,
$\left(\frac{2}{5}\right)=-1=\left(\frac{3}{5}\right)$ and
$d^2\equiv\pm1(5)$. Eventually, a contradiction is obtained by
applying the reciprocity law
$\left(\frac{5}{m}\right)=\left(\frac{m}{5}\right)$, if $5a+2$
and hence $m$ is odd, or
$\left(\frac{5}{n}\right)=\left(\frac{n}{5}\right)$, if $5a+3$ and
hence $n$ is odd.
\end{example}

\section{Moduli spaces yield orientation preserving equivalences}
\label{OrPreservSect}

Let us start out with a fairly general discussion of
Hodge isometries
$\varphi:\widetilde H(X,\IZ)\stackrel{\sim}{\to}\widetilde
H(X',\IZ)$ and the question when they are
orientation preserving.

Firstly, since $\varphi$ is an isomorphism of weight two Hodge
structures, we may choose generators $\sigma\in H^{2,0}(X)$ and
$\sigma'\in H^{2,0}(X')$ such that $\varphi(\sigma)=\sigma'$. The
oriented plane $\langle 1-\omega^2/2,\omega\rangle$,
where $\omega$ is a K{\"a}hler class, is completely
encoded by the complex line spanned by $\exp(i\omega)$. Moreover,
as $\varphi$ is an isometry, the image $\varphi(\exp(i\omega))$
is orthogonal to $\sigma'=\varphi(\sigma)$ and, therefore, of the
form
$$\varphi(\exp(i\omega))=\lambda\cdot\exp(b+ia)$$ with
$\lambda\in\IC^*$ and $a,b\in H^{1,1}(X',\IR)$.

Let us first compute the scalar $\lambda$. It will be expressed
as a linear combination of the degree zero parts of certain
natural Mukai vectors. We shall use the following short hand
$r:=\varphi(0,0,1)_0$, $\chi:=\varphi(1,0,1)_0$, and
$\chi_H:=\varphi(0,\omega,-\omega^2/2)_0$. Here, we are
anticipating the moduli space situation. Indeed, if $x\in X$ is a
closed point and $H\subset X$ an ample divisor with fundamental
class $\omega$, then $v(k(x))=(0,0,1)$, $v(\ko_X)=(1,0,1)$, and
$v(\ko_H)=(0,\omega,-\omega^2/2)$.

\begin{lemma}
$\lambda=\chi-r\left(\frac{\omega^2}{2}+1\right)+i\left(\chi_H+r\frac{\omega^2}{2}\right)$.
\end{lemma}

\begin{proof}
Write
$$\exp(i\omega)=v(\ko_X)-\left(\frac{\omega^2}{2}+1\right)\cdot
v(k(x))+i\left(v(\ko_H)+\frac{\omega^2}{2}\cdot v(k(x))\right)$$
and use $\lambda=\varphi(\exp(i\omega))_0$.
\end{proof}

Let us introduce the basic classes
$$u_0:=-r\cdot[\ko_X]+\chi\cdot[k(x)]~~{\rm
and}~~u_1:=-r\cdot[\ko_H]+\chi_H\cdot[k(x)]$$ as elements in the
Grothendieck group $K(X)$ (cf.\ \cite[Ex.\ 8.1.8]{HL}). Their
Mukai vectors are given by
$$v(u_0)=(-r,0,-r+\chi)~~{\rm
resp.}~~v(u_1)=(0,-r\omega,r\frac{\omega^2}{2}+\chi_H).$$

\begin{lemma}
\label{descrofaLemma}Suppose $r\ne0$. Then
$$|\lambda|^2\cdot a=
\varphi\left(\left(\frac{\chi_H}{r}+\frac{\omega^2}{2}\right)\cdot
v(u_0)
+\left(\left(\frac{\omega^2}{2}+1\right)-\frac{\chi}{r}\right)\cdot
v(u_1)\right)_2.$$
\end{lemma}

\begin{proof}
The assertion follows from \begin{eqnarray*}
|\lambda|^2\cdot a&=&{\rm
Im}(\bar\lambda\cdot\varphi(\exp(i\omega)))_2\\
&=&\left(\chi-r\left(\frac{\omega^2}{2}+1\right)\right)\varphi(\omega)_2-
\left(\chi_H+r\frac{\omega^2}{2}\right)\varphi\left(1,0,-\frac{\omega^2}{2}\right)_2
\end{eqnarray*}
\end{proof}

If $r=0$ and $\chi_H\ne0$, then $a$ cannot be written as the image
of a linear combination of $v(u_0)$ and $v(u_1)$.

\begin{proposition}\label{GeneralCrit}
Let $\varphi:\widetilde H(X,\IZ)\stackrel{\sim}{\to}\widetilde
H(X',\IZ)$ be a Hodge isometry with $\varphi(0,0,1)_0\ne0$.
Suppose $H\subset X$ is an ample divisor. Then $\varphi$ is
orientation preserving if and only if
$$\left(\frac{\chi_H}{r}+\frac{\omega^2}{2}\right)\cdot\varphi
(v(u_0))_2
+\left(\left(\frac{\omega^2}{2}+1\right)-\frac{\chi}{r}\right)\cdot
\varphi(v(u_1))_2$$ is contained in the positive cone
$\kc_{X'}\subset H^{1,1}(X',\IR)$.
\end{proposition}
Recall that the positive cone $\kc_{X'}$ is the connected
component of the cone of all classes $\alpha\in H^{1,1}(X,\IR)$
with $\alpha^2>0$ that contains the ample (respectively
K{\"a}hler) classes.
\begin{proof}
Using the notation introduced before, this linear combination of
$\varphi(v(u_0))_2$ and $\varphi(v(u_1))_2$ is (up to the positive
real scalar $| \lambda|$) the class $a$.

Now use the following easy facts:

i) Real and imaginary part of $\exp(ia)$ induce the same
orientation of the two positive directions of $\widetilde
H^{1,1}(X',\IR)$ as the natural one given by $\exp(i\omega')$ with
$\omega'$ an ample class if and only if $a\in\kc_{X'}$.

ii) Real and imaginary part of $\exp(ia)$ induce the same
orientation of the two positive directions in $\widetilde
H^{1,1}(X',\IR)$ as $\lambda\cdot \exp(b+ia)$ for any
$\lambda\in\IC^*$ and any $b\in H^{1,1}(X',\IR)$.
\end{proof}

\begin{remark}\label{twistOrPres}
Either by applying the above proposition or by any other means, it
is easy to check that the following standard equivalences are
orientation preserving: i) $F\mapsto f_*(F)$, where $f$ is an
isomorphism, ii) Line bundle twists $F\mapsto L\otimes F$ with
$L\in{\rm Pic}(X)$, iii) Shift functor $F\mapsto F[i]$, iv) Twist
functor $F\mapsto T_E(F)$, where $E\in\Db(X)$ is an arbitrary
spherical object (e.g.\ $\ko$).

All known (auto)equivalences can be written as compositions of the
above ones and the ones induced by universal families of stable
sheaves. Thus, in order to decide whether at least all known
(auto)equivalences are orientation preserving, it suffices to
consider the case of a fine moduli space and the equivalence it
induces.
\end{remark}

Let  us now consider a fine moduli space $X'=M(v)$ of stable
sheaves $E$ with $v(E)=v$, where $v\in\widetilde H^{1,1}(X,\IZ)$
with $\langle v,v\rangle=0$. Here stability is meant with respect
to a chosen ample divisor $H\subset X$.

Mukai has shown that in such a situation the moduli space $M(v)$
is, if not empty, again a K3 surface and that the Fourier--Mukai
transform $\Phi:\Db(X)\to\Db(M(v))$ with kernel the universal
sheaf $\ke=\ke(v)$ on $X\times M(v)$ is an equivalence (see
\cite{Mu}).

\begin{proposition}\label{proporpreserv}
The induced Hodge isometry
$$\xymatrix{\Phi_*:\widetilde H(X,\IZ)\ar[r]^-\sim&\widetilde
H(M(v),\IZ)}$$ is orientation preserving.
\end{proposition}

\begin{proof} The idea of the proof is of course to apply
Proposition \ref{GeneralCrit} and to show that the class $a$
introduced in the general context is contained in the positive
cone. We will however show a more precise result for the case of
a fine moduli space of stable sheaves of positive rank, namely
that in this case $a$ is in fact ample.

 Suppose $m\in M(v)$ corresponds to a stable sheaf $E$ on
$X$. Then $v=v(E)=(\rk(E),{\rm c}_1(E),\chi(E)-\rk(E))$, which we
will write as $(r,\ell,s)$.

The invariants $r$, $\chi$, and $\chi_H$ introduced above, can in
this special situation geometrically be interpreted as:
\begin{eqnarray*}
r=\Phi_*(v(k(x))_0=\rk(\ke|_{\{x\}\times
M(v)})=\rk(E)\\
\chi
=\Phi_*(v(\ko_X))_0=\rk(Rp_*\ke)=\chi(E),~~{\rm  and}\\
\chi_H=\Phi_*(v(\ko_H))_0=\rk(Rp_*(\ke|_{\{H\}\times
M(v)}))=\chi(E|_H). \end{eqnarray*}
 Now observe that that
twisting $E\mapsto E(nH)$ defines an isomorphism
$$\xymatrix{t(n):M(v)\ar[r]^-\sim& M(\exp(n[H])\cdot v)}$$ under
which the universal families can be compared by
$$({\rm id}_X\times t(n))^*\ke(\exp(n[H])\cdot v)\cong
q^*\ko(nH)\otimes\ke(v).$$ Clearly, $\Phi$ induces an orientation
preserving Hodge isometry if and only if $t(n)_*\circ\Phi\circ
(\ko(nH)\otimes(~~))$ does (cf.\ Remark \ref{twistOrPres}), but
the latter is nothing but the Fourier--Mukai transform with
respect to the universal family $\ke(\exp(n[H])\cdot v)$ on
$X\times M(\exp(n[H])\cdot v)$. (Similarly, the ampleness of the
class $a$ for $\ke$ is equivalent to the ampleness of the
analogous class with respect to $\ke(\exp(n[H])\cdot v)$.)

Thus, in order to prove the assertion we may first twist by
$\ko(nH)$ for $n\gg0$. Under the additional assumptions  that
$r>0$ we may thus arrange things such that
$\left(\frac{\chi_H}{r}+\frac{[H]^2}{2}\right)>0$ and
$\left(\left(\frac{[H]^2}{2}+1\right)-\frac{\chi}{r}\right)>0$.
Moreover, we may assume from the very beginning that all sheaves
$E\in M(v)$ are so positive that the standard GIT construction of
the moduli space applies directly. In particular, we may assume
that the hypothesis of Theorem 8.1.11 in \cite{HL} are fulfilled.
Again, one has to assume that $r>0$ .

Thus, the line bundle $$\kl_0:=\det(\Phi(u_0))\in{\rm
Pic}(M(v))$$ is ample. Its Chern class is ${\rm
c}_1(\kl_0)=\Phi_*(v(u_0))_2$.  (The line bundle $\kl_0$ is the
descent of the standard polarization of the Quot-scheme.)
Combined with \cite[Prop.\ 8.1.10]{HL} we deduce from this that
also the line bundle
$$\kl_1:=\det(\Phi(u_1))$$ is nef (see also the comments
in \cite[Sect.\ 8.2]{HL}). Its Chern class is ${\rm
c}_1(\kl_1)=\Phi_*(v(u_1))_2$. Hence,
$\left(\frac{\chi_H}{r}+\frac{\omega^2}{2}\right)\cdot\Phi_*
(v(u_0))_2
+\left(\left(\frac{\omega^2}{2}+1\right)-\frac{\chi}{r}\right)\cdot
\Phi_*(v(u_1))_2$ is  an ample class and thus contained in the
positive cone $\kc_{M(v)}$. Proposition \ref{GeneralCrit} yields
the assertion.

It remains to prove the assertion for $r=0$. Here, it will only be
shown that $a$ is contained in the positive cone (and not that it
is ample). A geometric proof can be given by viewing the moduli
space $M(v)$ as a relative moduli space over the linear system of
curves that occur as support of a stable sheaf (the case of
sheaves concentrated in points being trivial), but a more elegant
argument was suggested to us by Yoshioka. Instead of working with
the line bundles $\kl_i$, $i=0,1$, considered above he suggested
to work directly with the ample line bundle on the moduli space
that naturally occurs in the construction of Simpson. Here are
the details of the argument.

We know that $a\in\pm\kc_M$. In order to show that $a\in\kc_M$ it
suffices to show that $\langle a,\beta\rangle>0$ for one ample
class $\beta$ on $M$. By the very construction of the moduli
space {\`a} la Simpson, the line bundle $\kl:=\det(p_*(\ke\otimes
q^*u))$ is ample for $m\gg0$, where
$u:=\chi\cdot[\ko(mH)]-\chi(E(m))\cdot[\ko]$. Again, we suppose
from the very beginning that the sheaves have been twisted such
that the moduli space can directly be described as a GIT-quotient
of a certain Quot-scheme parametrizing all stable sheaves $[E]\in
M$.

Revisiting the proof of Lemma \ref{descrofaLemma} we find that the
class $a$ is a positive multiple of
$\left(\chi-r\left(\frac{\omega^2}{2}+1\right)\right)\varphi(\omega)_2-
\left(\chi_H+r\frac{\omega^2}{2}\right)\varphi\left(1,0,-\frac{\omega^2}{2}\right)_2$.
So if we introduce
$\gamma:=\left(\chi-r\left(\frac{\omega^2}{2}+1\right)\right)\omega-
\left(\chi_H+r\frac{\omega^2}{2}\right)\left(1,0,-\frac{\omega^2}{2}\right)
$, then $a=\varphi(\gamma)_2$ up to a positive scalar. One checks
that $\varphi(\gamma)_0=0$ by using
$\chi_H=-\langle(0,\omega,-\omega^2/2),v\rangle$. Similarly,
$\varphi(u)_0=0$. Hence,
$\langle\varphi(\gamma)_2,\varphi(u)_2\rangle=\langle\varphi(\gamma),\varphi(u)\rangle=\langle\gamma,
u\rangle$. A straightforward calculation shows that
$\langle\gamma,u\rangle>0$ which concludes the proof.
\end{proof}

\begin{remark}
  The case of a fine moduli space of stable sheaf of positive rank
  is special in the sense that the
class $a$ given by $\Phi_\ke^H(\exp(i[H]))=\lambda\cdot\exp(b+ia)$
is not only contained in the positive cone, but is in fact ample.
This does not hold true for arbitrary Fourier--Mukai equivalences
as is shown e.g.\ by the twist with respect to the spherical
object $\ko_C(k)$, where $C\subset X$ is a smooth rational curve.
\end{remark}

\begin{remark}
We will use the main result of this section not only in the case
of fine moduli spaces, but as well for coarse moduli spaces. More
precisely, the $(1,\alpha_B)$-twisted universal family $\ke$ on
$X\times M(v)$ induces an  orientation preserving Hodge isometry
$\Phi_{*}^{0,B}:\widetilde H(X,\IZ)\stackrel{\sim}{\to}\widetilde
H(M(v),B,\IZ)$.

The above proof carries over to coarse moduli spaces by the
following rather standard arguments. The main idea is to use the
fact that the line bundles $\kl_i$ and $\kl$ are also defined for
coarse moduli spaces (even semi-stable sheaves are allowed).
Attention has to be paid to problems with the twist on the level
of cohomology.

Here are some of the details. The moduli space $M(v)$ is
constructed as a quotient $\pi:R^s\to M(v)$ of a certain open
subset $R^s$ of a certain Quot-scheme. In particular, there
always exists a universal sheaf $\tilde \ke$ on $X\times R^s$. In
fact, even in the case of a fine moduli space one first
constructs $\kl_0$ on $R^c$. Then one shows that it  is the
restriction of an ample line bundle on the ambient Quot-scheme
and that its restriction to $R^s$ descends to an ample line
bundle on $M(v)$.

If a twisted universal sheaf $\ke$ is given by $\{X\times U_i,
\varphi_{ij}\in\ko^*(X\times U_{ij})\}$ then $\tilde\ke$ can be
thought of as given by
$\{X\times\pi^{-1}(U_i),\tilde\varphi_{ij}:=(1\times\pi)^*\varphi_{ij}\cdot
\lambda_{ij}\}$. Here, $\lambda_{ij}\in\ko^*(X\times
\pi^{-1}(U_{ij}))$ satisfy
$\lambda_{ij}\cdot\lambda_{jk}\cdot\lambda_{ki}=(1\times\pi)^*\alpha_{ijk}^{-1}$.

Recall that also the sheaf  $\ke_B$ needed to define ${\rm
ch}^B(\ke)$ was obtained by untwisting with invertible functions
$\exp(a_{ij})$ satisfying a similar cocycle condition as the
$\lambda_{ij}$ but already on $X\times M(v)$. The functions
$a_{ij}$ are defined  on $X\times U_{ij}$, but they are not
holomorphic. Thus, one finds that $\tilde\ke$ and
$(1\times\pi)^*\ke_B$ differ by a (non-holomorphic) line bundle
(the one given by the transition functions
$\lambda_{ij}\cdot\exp(-a_{ij})$). For the same reason that shows
that for fine mo\-duli spaces the line bundles $\kl_i$ do not
depend on the chosen universal family (which is only unique up to
twist by line bundles on the moduli space) one concludes that
$\pi^*\Phi_*^{0,B}(u_i)_2=\pi^*{\rm c}_1(\kl_i)$ as well as
$\pi^*\Phi_*^{0,B}(u)_2=\pi^*{\rm c}_1(\kl)$ (using the notation
of the proof of Proposition \ref{proporpreserv}). This is enough
to conclude also in the case of coarse moduli spaces.
\end{remark}

\section{(Twisted) derived
isometries}\label{TwistedisometriesSection}

Let us begin with a few comments on the group
$\OO(\widetilde\Lambda)$ of isometries of the lattice
$\widetilde\Lambda=\Lambda\oplus U$. Classical results due to C.\
T.\ C.\ Wall show that the three natural subgroups
$\OO(\Lambda),\OO(U),\exp(\Lambda)\subset\OO(\widetilde \Lambda)$
generate $\OO(\widetilde\Lambda)$ (see \cite{W}). Here, $\exp(B)$
with $B\in\Lambda$ acts by multiplication with $1+B+B^2/2$ in
$\widetilde\Lambda=H^*(M,\IZ)$.  Clearly, the subgroups $\OO(U)$
and $\OO(\Lambda)$ commute. For $g\in\OO(\Lambda)$ and
$\exp(B)\in\exp(\Lambda)$ one has
$$g\circ\exp(B)=\exp(g(B))\circ g.$$
Hence, any element $g\in\OO(\widetilde\Lambda)$ can be written as
$$g=g_1\circ g_2~~{\rm
with}~~g_1\in\langle\OO(U),\exp(\Lambda)\rangle~~{\rm and
}~~g_2\in\OO(\Lambda).$$ Also note that $j=-{\rm id}_U$ commutes
with $\OO(U),\OO(\Lambda)$ and satisfies $\exp(B)\circ
j=j\circ\exp(-B)$. The group of orientation preserving isometries
$\OO_+(\widetilde\Lambda)$ can similarly be generated by
$\exp(\Lambda)$, $\OO_+(\Lambda)$, and $\OO_+(U)$.

 In the following, we will be interested in the
following subsets
$$\begin{array}{rcl}
H&:=&\{g\in\OO_+(\widetilde\Lambda)~|~g(Q)\cap Q\ne\emptyset\}\\
H_{\rm alg}&:=&\{g\in\OO_+(\widetilde\Lambda)~|~g(Q_{\rm alg})\cap
Q_{\rm alg}\ne\emptyset\},
\end{array}$$
where $Q\subset\IP(\Lambda_\IC)$ is the K3 surface period domain
$\{x~|~x^2=0,(x,\bar x)>0\}$ and $Q_{\rm alg}\subset Q$ is the
dense subset of periods of algebraic K3 surfaces. In other words,
$Q_{\rm alg}$ is the set of those $x\in Q$ for which there exists
a class $B\in x^\perp\cap\Lambda$ with $B^2>0$.

\begin{proposition}\label{untwistedisom}
Both sets $H_{\rm alg}\subset H\subset\OO_+(\widetilde\Lambda)$
contain the generating subgroups
$\OO_+(\Lambda),\OO_+(U),\exp(\Lambda)\subset\OO_+(\widetilde\Lambda)$
and thus generate $\OO_+(\widetilde\Lambda)$. However,
$H\ne\OO_+(\widetilde\Lambda)$ or, equivalently, the subsets
$H_{\rm alg}$ and $H$ do not form subgroups of $\OO_+(\widetilde
\Lambda)$.
\end{proposition}

\begin{proof}
In fact, $\OO_+(\Lambda)$ and $\OO_+(U)$ both respect the period
domain $Q$. If $x\in Q_{\rm alg}\cap B^\perp$ for a given
$B\in\Lambda$, then $\exp(B)x=x$ and hence $\exp(B)\in H_{\rm
alg}$. The existence of $x$ follows from the fact that $B^\perp$
contains a positive plane.

In order to show the second assertion, it suffices to construct
one $g\in\OO_+(\widetilde\Lambda)$ with $g(Q)\cap Q=\emptyset$.

Choose $B_0,B_1\in\Lambda$ such that $\langle
B_0,B_1\rangle\subset \Lambda_\IR$ is a positive plane. Then
consider the isometry $g:=\exp(B_1)\circ i\circ\exp(B_0)$ where
$i\in\OO_+(U)$ is the isometry that maps the generators $u_2$ and
$u_1$ of $H^0(M,\IZ)$ respectively $H^4(M,\IZ)$ to $-u_1$
respectively $-u_2$.

Suppose $g(x)\in Q$ for some $x\in Q$. Writing $g(x)=
-(B_0,x)u_2+(x-(B_0,x)B_1)+(B_1-\frac{B_1^2}{2}B_0,x)u_1$ then
shows that $x\in B_0^\perp\cap B_1^\perp$. This is impossible, as
real and imaginary part of $x$ together with $B_0$ and $B_1$ would
span a positive four-space in $\Lambda_\IR$ that does not exist.
\end{proof}


Let us now introduce the following notation:
$$\begin{array}{rcl}
H'&:=&\{g\in\OO_+(\widetilde\Lambda)~|~g(\exp(\Lambda_\IQ) \cdot Q)\cap (\exp(\Lambda_\IQ)\cdot Q)\ne\emptyset\}\\
H'_{\rm
alg}&:=&\{g\in\OO_+(\widetilde\Lambda)~|~g(\exp(\Lambda_\IQ)\cdot
Q_{\rm alg})\cap (\exp(\Lambda_\IQ)\cdot Q_{\rm
alg})\ne\emptyset\}.
\end{array}$$

\begin{proposition}\label{twistedIso}
$H'_{\rm alg}=H'=\OO_+(\widetilde\Lambda)$. In particular, both
sets are groups.
\end{proposition}

\begin{proof}
We will actually show slightly more, namely that for any
$g\in\OO(\widetilde\Lambda)$ one has $g(Q_{\rm
alg})\cap(\exp(\Lambda_\IQ)\cdot Q_{\rm alg})\ne\emptyset$.

Fix a basis $x_1,\ldots,x_{22}\in\Lambda$. Then any
$g\in\OO(\widetilde\Lambda)$ can be written as
$$g(x)=\sum_{i=1}^{22}\lambda_ix_i+\mu_1u_1+\mu_2u_2$$
with $\lambda_i,\mu_j\in\IZ$ depending linearly on $x\in\Lambda$.
(Here, $u_2$ and $u_1$ span $H^0(M,\IZ)$ respectively
$H^4(M,\IZ)$.) In other words, there exist elements
$A_i,B_1,B_2\in\Lambda$ with
$$g(x)=\sum_{i=1}^{22}(A_i,x)x_i+(B_1,x)u_1+(B_2,x)u_2$$
for any $x\in\Lambda$. Of course, by linear extension the same
formula holds for any $x\in\Lambda_\IC$.

Now choose
a positive plane inside $B_2^\perp\subset\Lambda_\IQ$ and an
orthogonal basis $x_1,x_2\in \Lambda_\IQ$ of it with $x_1^2=x_2^2$.
Then set $x:=x_1+ix_2$. With this definition $x\in Q$ corresponds to
a K3 surface of Picard number $20$ which is necessarily projective,
i.e.\ $x\in Q_{\rm alg}$.

As $(B_1,x)=(B_1,x_1)+i(B_1,x_2)\in\IQ(i)$, one finds
$B\in\Lambda_\IQ$ such that $(B_1,x)=(B,y)$, where $y\in Q$ is
given by $g(x)=y+(B_1,x)u_1$. As $x\in Q_{\rm alg}$ immediately
yields $y\in Q_{\rm alg}$, this proves
$g(y)\in\exp(\Lambda_\IQ)\cdot Q_{\rm alg}$.

This then shows that $g(Q_{\rm alg})\cap(\exp(\Lambda_\IQ)\cdot
Q_{\rm alg})\ne\emptyset$.
\end{proof}

\begin{remark}\label{twistedIsoRem}
i) Note that the K3 surface corresponding to the period $x$
considered in the above proof, which is defined over $\IQ$, has
maximal Picard number $\rho=20$.

ii) The above result can be improved to the following statement: For
any $g\in \OO_+(\widetilde\Lambda)$ one has $g(Q_{\rm
alg})\cap(\exp(\Lambda)\cdot Q_{\rm alg})\ne\emptyset$. Indeed, we
may choose integral classes $x_1,x_2\in B_2^\perp$ with
$x_1^2=x_2^2=2$ and $\langle x_1,x_2\rangle=0$. (Note that, after
applying isometries of $\widetilde\Lambda$, any primitive vector can
be assumed to be contained in one copy $U$ in the decomposition
$\widetilde \Lambda=U^{\oplus 4}\oplus(-E_8)^{\oplus 2}$. Thus
$B_2^\perp$ contains in particular two other copies of the
hyperbolic plane, which ensures the existence of $x_1,x_2$.) Then,
real and imaginary parts $y_1,y_2$ of $y$, defined as in the proof,
have the same properties. It is now possible to find an integral
element $B\in\widetilde\Lambda$ with $\langle B,y_i\rangle=\langle
B_1,x_i\rangle$. Indeed, the primitive sublattice generated by
$y_1,y_2$ is isometric to the sublattice generated by $e_1+e_2$ and
$e_1'+e_2'$ inside the direct sum of two copies of the hyperbolic
plane $U\oplus U'\subset\widetilde\Lambda$, which form a direct
summand of $\widetilde\Lambda$ . Standard results of Nikulin show
that this isometry can be extended to an isometry of
$\widetilde\Lambda$. The inverse image under this isometry of
$(\langle B_1,x_1\rangle e_1)+(\langle B_1,x_2\rangle e'_1)$ can be
taken for $B$. Note that $x=x_1+ix_2$ is again algebraic.
\end{remark}

We continue to present a K3 surface by a complex structure $I$ on
the fixed manifold $M$. In particular, a marking $H^2(X,\IZ)\cong
H^2(M,\IZ)\cong\Lambda$ is automatically given.

\begin{definition}
An element $g\in\OO_+(\widetilde\Lambda)$ is called a
\emph{derived isometry} if there exist two algebraic K3 surfaces
$X$ and $X'$ and a Fourier--Mukai equivalence
$\Phi:\Db(X)\cong\Db(X')$ with $\Phi_*=g$.

An element $g\in\OO_+(\widetilde\Lambda)$ is called a
\emph{twisted derived isometry} if there exist two algebraic K3
surfaces $X$, $X'$, B-fields $B\in H^2(X,\IZ)$, $B'\in
H^2(X',\IZ)$ and a twisted Fourier--Mukai equivalence
$\Phi:\Db(X,\alpha_B)\cong \Db(X',\alpha_{B'})$ such that
$\Phi_*^{B,B'}=g$.
\end{definition}

In the twisted as well as in the untwisted situation, the
difficult question seems to be whether any (twisted) derived
equivalence induces a (twisted) derived isometry in this sense,
i.e.\ whether it is orientation preserving. For all known examples
this is the case, as was shown in the previous section.

In order to emphasize the difference between these two notions we
sometimes speak of untwisted derived isometry in the first case.

\begin{corollary}
The set $G\subset\OO_+(\widetilde \Lambda)$ of all derived
isometries generates $\OO_+(\widetilde\Lambda)$. However,  $G$ is
not a group, i.e.\ $G\ne \OO_+(\widetilde\Lambda)$.
\end{corollary}

\begin{proof}
By Borcea's result any element in $\OO_+(\Lambda)$ is realized as
an isomorphism of K3 surfaces. Hence $\OO_+(\Lambda)\subset G$.

The subgroup $\OO_+(U)=\OO_+((H^0\oplus H^4)(M,\IZ))$ consists of
the two elements ${\rm id}$ and $i$, where $i$ is as in the proof
of Proposition \ref{untwistedisom}. If $X$ is an arbitrary K3
surface, then the reflection $T_\ko$ associated to the spherical
object $\ko$ acts as $i$ on cohomology. Hence, $G$ contains
$\OO_+(U)$.

 Eventually, any $\exp(B)$ with $B\in \Lambda$ is
contained in $G$. One way to see this is to choose an algebraic
K3 surface $X$ for which $B$ is of type $(1,1)$, i.e.\ $B={\rm
c}_1(L)$ for some $L\in{\rm Pic}(X)$. Then the Fourier--Mukai
equivalence given by $\otimes L$ acts as multiplication with
$\exp(B)={\rm ch}(L)$ on cohomology.

Thus, $G$ contains all three subgroups $\OO_+(\Lambda)$,
$\OO_+(U)$, and $\{\exp(B)\}_{B\in\Lambda}$ and hence generates
$\OO_+(\widetilde \Lambda)$.

The last assertion follows from Proposition \ref{untwistedisom}
and the fact that $G\subset H_{\rm alg}$.
\end{proof}

 The positive result in the
twisted case we are going to present confirms once more the
difference between the twisted and the untwisted situation. It is
not a direct consequence of Proposition \ref{twistedIso}, as it
uses C\u{a}ld\u{a}raru's conjecture for large Picard number which
will be established only in the next section.

\begin{proposition}\label{twistedisoformgroup}
Let $G'\subset\OO_+(\widetilde \Lambda)$ be the set of all twisted
isometries. Then $G'=\OO_+(\widetilde\Lambda)$.
\end{proposition}

\begin{proof}
Let $g\in \OO(\widetilde\Lambda)$. By Proposition
\ref{twistedIso} there exist twisted algebraic K3 surfaces
$(X,\alpha_B)$ and $(X',\alpha_{B'})$ such that $g$ defines an
Hodge isometry $\widetilde H(X,B,\IZ)\cong \widetilde
H(X',B',\IZ)$. If one could use  C\u{a}ld\u{a}raru's conjecture
as stated in \ref{Calconjref}, then one immediately would deduce
the existence of a twisted derived Fourier--Mukai equivalence
$\Db(X,\alpha_B)\cong\Db(X',\alpha_{B'})$ with $\Phi_*^{B,B'}= g$.

Fortunately, in the proof of Proposition \ref{twistedIso} we
actually constructed $X$ and $X'$ of Picard number $\rho\geq12$
(cf.\ Remark \ref{twistedIsoRem}) and for those a variant of
C\u{a}ld\u{a}raru's conjecture, that suffices for this argument,
will be proved in the next section.

One could also argue without evoking C\u{a}ld\u{a}raru's
conjecture at all by using ii), Remark \ref{twistedIsoRem}.
Indeed, this remark shows that for any $g\in \OO(\widetilde\Lambda)$
there exist two algebraic K3 surfaces $X$ and $X'$, an equivalence
$\Phi:{\rm D}^{\rm b}(X)\cong{\rm D}^{\rm b}(X')$, and
an integral B-field $B'\in H^2(X',\IZ)$ such that
$g=\exp(B')\cdot\Phi_*$.
\end{proof}

\section{C\u{a}ld\u{a}raru's conjecture for large Picard
number}\label{SectionlargePicard}

Let us start out by recalling the following

\begin{proposition}\label{MukaiProp}{\bf (Mukai)} Suppose $X_1$ and $X_2$ are two
K3 surfaces with Picard number $\rho(X_i)\geq12$. Then up to sign
any Hodge isometry $T(X_1)\cong T(X_2)$ is induced by an
isomorphism $X_1\cong X_2$.
\end{proposition}

\begin{remark}
i) As $(H^0\oplus H^4)(M,\IZ)$ forms a hyperbolic plane contained
in the orthogonal complement of any transcendental lattice
$T(X)\subset \widetilde H(M,\IZ)$, any Hodge isometry $T(X_1)\cong
T(X_2)$ can be extended to a Hodge isometry $\widetilde
H(X_1,\IZ)\cong\widetilde H(X_2,\IZ)$. (This follows from
\cite[Thm.1.14.4]{Nikulin}, see Remark \ref{NikulinRem},
and no additional assumption on the
Picard number is needed for this.)

It is this fact that allows one to phrase Orlov's result in terms
of the transcendental lattices as: $\Db(X_1)\cong\Db(X_2)$ if and
only if $T(X_1)\cong T(X_2)$. Thus, as isomorphic K3 surfaces
have equivalent derived categories, this result of Mukai's proves
in particular Orlov's result under the additional assumption
$\rho\geq12$.

ii) Using the Global Torelli theorem, the proof of the above
proposition reduces to the purely lattice theoretical problem to
extend a given Hodge isometry $T(X_1)\cong T(X_2)$ to a Hodge
isometry $H^2(X_1,\IZ)\cong H^2(X_2,\IZ)$. This is possible by
results of Nikulin \cite{Nikulin}, as the orthogonal complement of
$T(X_i)$ in $H^2(X_i,\IZ)$, i.e.\ the Picard group, is big enough
by assumption.

Note that the occurrence of a possible sign is missing in Prop.\
6.3 in \cite{Mu}.
\end{remark}

 The principal result of this section establishes
C\u{a}ld\u{a}raru's conjecture under a similar assumption on the
Picard number. The slight difference to Mukai's original result
is that in the twisted case only derived equivalence holds and
not isomorphism of the twisted K3 surfaces. For examples see
Section \ref{FMcounting}.

Before coming to this, let us show how to reduce the twisted
problem for large Picard number to the untwisted situation.

\begin{proposition}\label{twistedisuntwisted}
Let $(X,\alpha_B)$ be a  twisted K3 surface with $\rho(X)\geq12$.
Then there exists a K3 surface $Z$ and a Fourier--Mukai equivalence
$\Phi:\Db(Z)\cong\Db(X,\alpha_B)$. Moreover, $\Phi$ can be chosen
such that $\Phi^{0,B}_*$ is orientation preserving.
\end{proposition}

\begin{proof}
The proof follows almost directly from results of Mukai and
C\u{a}ld\u{a}raru, so we will be brief.

As $\rk(T(X,B))=\rk(T(X))\leq10$, \cite[Thm.1.14.4]{Nikulin}
applies  (see Remark \ref{NikulinRem})
and shows that there exists a primitive embedding
$T(X,B)\hookrightarrow\Lambda$ (see the discussion in Remark
\ref{NikulinRem}). By the surjectivity of the period map, one
finds a K3 surface $Z$ and a Hodge isometry $T(Z)\cong T(X,B)$.

The composition of this Hodge isometry with $\exp(-B)$ yields an
embedding $i:T(Z)\hookrightarrow T(X)$ which is the kernel of the
natural map $\alpha:T(X)\to\IZ/n\IZ$ defined by $\alpha=\alpha_B$,
where $n$ is the order of $\alpha\in{\rm Br}(X)$. Choose an
element $\ell\in T(X)$ with $\alpha(\ell)=1\in\IZ/n\IZ$ and
denote by $t\in T(Z)$ the element with $ i(t)=n\cdot\ell$.

Next, we use arguments of Mukai (see Section 6 in \cite{Mu}). He
shows that there exists a compact, smooth, two-dimensional moduli
space $M$ of stable sheaves on $Z$ such that the inclusion
$\varphi:T(Z)\to T(M)$, which is defined in terms of a
quasi-universal family, maps $t$ to an element in $T(M)$ which is
divisible by $n$ and such that $(1/n)\varphi(t)$ generates the
quotient ${\rm Coker}(\varphi)\cong\IZ/n\IZ$.

C\u{a}ld\u{a}raru continued Mukai's discussion and showed that the
Brauer class $\beta\in{\rm Br}(M)$ defined by ${\rm
Coker}(\varphi)\cong\IZ/n\IZ$, $(1/n)\varphi(t)\mapsto 1$ is the
obstruction class for the existence of a universal sheaf.
Moreover, he showed that a $(1,\beta^{-1})$-twisted universal
sheaf $\ke$ exists and that the induced Fourier--Mukai functor
defines an equivalence $\Phi:\Db(Z)\cong\Db(M,\beta^{-1})$
which is orientation preserving (see Prop.\ \ref{OrPreservSect}).
C\u{a}ld\u{a}raru checked this equivalence by applying the
standard criterion due to Bondal and Orlov \cite{BO,Bridgeland}
testing a Fourier--Mukai functor on points.

In the penultimate step, one remarks that the two inclusions
$i:T(Z)\hookrightarrow T(X)$ and $\varphi:T(Z)\hookrightarrow
T(M)$ can be identified via an isomorphism $\psi:T(X)\cong T(M)$
which sends $\ell=(1/n)i(t)$ to $(1/n)\varphi(t)$. This yields a
commutative diagram
$$\xymatrix{T(Z)\ar@{=}[d]\ar@{^(->}[r]^i&T(X)\ar[d]^\psi\ar@{>>}[r]&\IZ/n\IZ\ar@{=}[d]\\
T(Z)\ar@{^(->}[r]^\varphi&T(M)\ar@{>>}[r]&\IZ/n\IZ.}$$

Eventually, we use Mukai's result Proposition \ref{MukaiProp} to
find an isomorphism $f:M\to X$ with $f^*|_{T(X)}=\pm\psi$.
Let us first suppose that  $f^*|_{T(X)}=-\psi$.
In
view of the commutativity of the above diagram this morphism
satisfies  $f^*\alpha=\beta^{-1}$. Which
shows that there exists a Fourier--Mukai equivalence
$$\xymatrix{\Db(Z)\ar[r]^-\sim_-{\Phi}&\Db(M,\beta^{-1})\ar[r]^-\sim_{f_*}&\Db(X,\alpha)}$$
that preserves the orientation. If $f^*|_{T(X)}=\psi$, then use
the composition
$$
\xymatrix{\Db(Z)\ar[r]^-\sim_-{\Psi}&\Db(M,\beta^{-1})\ar[r]^-\sim_{f_*}&\Db(X,\alpha).}$$
Here, $\Psi$ is the Fourier--Mukai transform with kernel $\ke^*$
 that is obtained by dualizing the
$(1,\beta^{-1})$-twisted universal family $\ke$ on $Z\times M$.
Thus, $\ke^*$ is a $(1,\beta)$-object. Using the standard
criterion it is easy to see that the Fourier--Mukai transform
$\Db(Z)\cong\Db(M,\beta)$ induced by it is also an equivalence.
This is essentially due to the fact that
$\Ext^i(\ke_x,\ke_y)\cong\Ext^{i}(\ke^*_y,\ke^*_x)$ for all
$x,y\in Z$ . (The statement reminds of the fact that
$\Phi_{\ke^*}:\Db(X)\cong\Db(Y)$ is an equivalence if
$\Phi_\ke:\Db(X)\cong\Db(Y)$ is one. A proof of this fact not
relying on the point criterion can be found in \cite{Or2}.)

A quick look at the induced cohomological Fourier--Mukai transform
reveals that $\Psi_{*}$ is as well orientation preserving. Indeed,
if $\Phi_*(\exp(i\omega))=\lambda\exp(b+ia)$ and
$\Psi_{*}(\exp(i\omega))=\lambda^*\exp(b^*+ia^*)$ then
$\lambda^*=\bar\lambda$ and $\lambda^*\exp(b^*+ia^*)=-\overline
{\lambda\exp(b+ia)}$ (we use the notation of Section
\ref{OrPreservSect}). Hence, the imaginary part of the degree two
part of the image of $\exp(i\omega)$ does not change.
\end{proof}

\begin{corollary}
Let $(X,\alpha)$ and $(X',\alpha')$ be  twisted K3 surfaces.
Assume $\rho(X)=\rho(X')\geq 12$. Then any equivalence
$\Db(X,\alpha)\cong\Db(X',\alpha')$ is of Fourier--Mukai type.
\end{corollary}

\begin{proof}
The proposition provides Fourier--Mukai equivalences
$\Phi:\Db(Z)\cong\Db(X,\alpha)$ and
$\Phi':\Db(Z')\cong\Db(X',\alpha')$. If
$\Psi:\Db(X,\alpha)\cong\Db(X',\alpha')$ is any equivalence, then
${\Phi'}^{-1}\circ\Psi\circ\Phi:\Db(Z)\cong\Db(Z')$ is an
equivalence between untwisted derived categories and hence, due
to the result of Orlov, of Fourier--Mukai type. This is enough to
conclude that also $\Psi$ is of Fourier--Mukai type.
\end{proof}

There is a natural Hodge isometry $\widetilde
H(X,B,\IZ)\cong\widetilde H(X,-B,\IZ)$ provided by $j=(-{\rm
id}_{H^0})\oplus {\rm id}_{H^2}\oplus (-{\rm id}_{H^4})$. However,
it is not preserving the orientation of the positive directions.
So, according to the modified C\u{a}ld\u{a}raru conjecture
\ref{Calconjref} we do not expect this Hodge isometry to lift to
a Fourier-Mukai equivalence
$\Db(X,\alpha_B)\cong\Db(X,\alpha_{-B})$. In fact, in general the
two twisted derived categories $\Db(X,\alpha)$ and
$\Db(X,\alpha^{-1})$ are probably not equivalent (except if
$\alpha$ is of order two). However, if the Picard group is big,
this is true and follows from the proof of the proposition:

\begin{corollary}\label{dualtwisted}
Let $(X,\alpha)$ be a twisted K3 surface with $\rho(X)\geq 12$.
Then there exists an orientation preserving
Fourier--Mukai equivalence
$\Db(X,\alpha)\cong\Db(X,\alpha^{-1})$.\qqed
\end{corollary}

\begin{remark}
The arguments given in the proof also show that any Fourier--Mukai
equivalence $\Phi:\Db(X,\alpha)\cong\Db(X',\alpha')$ induces a
Fourier--Mukai equivalence
$\Db(X,\alpha^{-1})\cong\Db(X',{\alpha'}^{-1})$.

We emphasize that the equivalence
$\Db(X,\alpha)\cong\Db(X,\alpha^{-1})$ is by no means canonical,
as it depends on the choice of the intermediate K3 surface $Z$
and the kernel yielding the equivalence
$\Db(Z)\cong\Db(X,\alpha)$. This becomes evident if one considers
a case where $\alpha$ is trivial, $X=Z$, and $\Db(Z)\cong\Db(X)$
is given by a kernel $\ke$ different from the diagonal. Then the
autoequivalence of $\Db(X)$ constructed by the above methods
would be $\Phi_\ke^2$, which is in no way natural.
\end{remark}

Let us now come to the main result of this section.

\begin{proposition}\label{MainSect6}
Let $X_1$, $X_2$ be two algebraic K3 surfaces with rational
B-fields $B_1$ respectively $B_2$. Assume $\rho(X_1)\geq 12$.

If $g:\widetilde H(X_1,B_1,\IZ)\cong \widetilde H(X_2,B_2,\IZ)$
is an orientation preserving Hodge isometry, then there exists a Fourier--Mukai
equivalence
$$\Phi:\Db(X_1,\alpha_{B_1})\cong\Db(X_2,\alpha_{B_2})$$ with
$\Phi^{B_1,B_2}_*=g$.
%
\end{proposition}

\begin{proof}
The previous proposition provides K3 surfaces $Z_1$ and $Z_2$
together with
Fourier--Mukai equivalences ${\Phi_1}:\Db(Z_1)\cong
\Db(X_1,\alpha_{B_1})$ and
${\Phi_2}:\Db(Z_2)\cong\Db(X_2,\alpha_{B_2})$
both preserving the orientation of the positive directions
in cohomology. On the level of
cohomology this yields an orientation preserving
Hodge isometry (cf.\ Proposition
\ref{HyieldsD})
$$\xymatrix{h:\widetilde H(Z_1,\IZ)\ar[r]^-{{\Phi^{0,B_1}_{1*}}}&\widetilde
H(X_1,B_1,\IZ)\ar[r]^g&\widetilde
H(X_2,B_2,\IZ)\ar[r]^-{(\Phi_{2*}^{0,B_2})^{-1}}&\widetilde
H(Z_2,\IZ).}$$
Hence, by results of
Mukai, Orlov et al, there exists an equivalence
$\Psi:\Db(Z_1)\cong\Db(Z_2)$ with $\Psi_*=h$.
\end{proof}

\begin{corollary}
Let $X_1$ and $X_2$ be algebraic K3 surfaces with
$\rho(X_i)\geq12$ and $B_i\in H^2(X_i,\IQ)$, $i=1,2$. Then the
following conditions are equivalent.

{\rm i)} There exists a Fourier--Mukai equivalence
$\Phi:\Db(X_1,\alpha_{B_1})\cong\Db(X_2,\alpha_{B_2})$.

{\rm ii)} There exists a Hodge isometry $\widetilde
H(X_1,B_1,\IZ)\cong\widetilde H(X_2,B_2,\IZ)$.

{\rm iii)} There exists a Hodge isometry $T(X_1,B_1)\cong
T(X_2,B_2)$.
\end{corollary}

\begin{proof}
By Proposition \ref{HyieldsD}, it is clear that ii) and iii)
follow from i) (without any assumption on the Picard number).
That ii) implies i) follows from Proposition \ref{MainSect6} and
Corollary \ref{dualtwisted}.

In order to see that iii) implies ii) one remarks that the
natural embedding
$$T(X,B)\hookrightarrow\widetilde
H(X,B,\IZ)\cong\widetilde\Lambda$$ is unique by the results of
Nikulin \cite[Thm.1.14.4]{Nikulin} (see Remark
\ref{NikulinRem}). Hence, any Hodge isometry
between the twisted transcendental lattices extends to a Hodge
isometry of the full twisted Hodge structures.
\end{proof}

The following result answers a question of C\u{a}ld\u{a}raru (see
5.5.3 in \cite{Cal})  affirmatively in the case of large Picard
number, although the answer in the general case should be negative
(see Example \ref{counter}). At the same time, it generalizes
Corollary \ref{dualtwisted}.

\begin{corollary}
Let $X$ be a K3 surface with $\rho(X)\geq12$ and $\alpha\in{\rm
Br}(X)$. If $k$ is prime to the order of $\alpha$, then there
exists a Fourier--Mukai equivalence
$\Db(X,\alpha)\cong\Db(X,\alpha^k)$.
\end{corollary}

\begin{proof}
Just observe that the kernel of
$\alpha:T(X)\twoheadrightarrow\IZ/n\IZ$ and
$\alpha^k:T(X)\twoheadrightarrow\IZ/n\IZ$ are Hodge isometric and
apply the previous corollary.
\end{proof}

\begin{remark}\label{Brauerorderdifferent}
A priori the order of the two Brauer classes $\alpha_{B_1}$ and
$\alpha_{B_2}$ in the proposition might be different (cf.\ Remark
\ref{notthesameorder}). Although we have seen earlier that Brauer
classes on the same K3 surface defining equivalent twisted
derived categories are of the same order (see Remark
\ref{ordertrans}).
%
\end{remark}

\begin{remark}\label{NikulinRem}
As it turns out, all results of this section hold true under  a
slightly technical but weaker lattice theoretical condition. Let
us briefly explain this here. First we recall the original result of
Nikulin  \cite[Thm.1.14.4]{Nikulin}
as given in \cite{Mo}: Let $T$ and $\Lambda$ be
an even respectively an even unimodular lattice such that
$t_\pm<s_\pm$ and
\begin{equation}\label{latticecond}
  \ell(T)\leq \rk(\Lambda)-\rk(T)-2.
\end{equation}
Then there exists a unique (up to isometries of
$\Lambda$) primitive embedding of $T$ into $\Lambda$.
Here, $(t_+,t_-)$ and $(s_+,s_-)$ are the signatures of
the lattices $T$ respectiveley $\Lambda$
and $\ell(T)$ is the minimal number of generators
of the discriminant group $A_T=T\dual/T$.

If $\rk(T)\leq 10$ then (\ref{latticecond}) holds true with
$\Lambda$ the K3 lattice. But clearly there are other
cases when these assumptions are satisfied.
\end{remark}

\section{Counting twisted Fourier--Mukai
partners}\label{FMcounting}

In the untwisted as well as in the twisted case it has been shown
that the number of non-isomorphic Fourier--Mukai partners of a
given algebraic (twisted) K3 surface is always finite. But the
number itself can be arbitrarily large. More precisely, it has
been shown in \cite{Og,St} that for any $N$ there exist pairwise
non-isomorphic K3 surfaces $X_1,\ldots,X_N$ with
$$\Db(X_1)\cong\ldots\cong\Db(X_N).$$

However, in the untwisted case the Picard number of these K3
surfaces has to be small, i.e.\ $\rho(X_i)<12$. Indeed, if
$\Db(X_1)\cong\Db(X_2)$ and $\rho(X_i)\geq 12$, then due to
Proposition \ref{MukaiProp} one automatically has $X_1\cong X_2$.
Moreover, following \cite{Oguisonew}, this also holds true if
$\rho\geq3$ and the determinant of ${\rm Pic}(X)$ is square free.
For the calculation of the number of Fourier--Mukai partners in the
untwisted case see \cite{Oguisonew,St}.

 Passing to the twisted situation allows
one to construct arbitrarily many (say $N$) pairwise
non-isomorphic twisted K3 surfaces
$(X_1,\alpha_1),\ldots,(X_N,\alpha_N)$ of large Picard number.

More precisely, one has

\begin{proposition}
For any $N$ there exist $N$ pairwise non-isomorphic algebraic K3
surfaces $X_1,\ldots,X_N$ of Picard number $\rho(X_i)=20$, endowed
with Brauer classes $\alpha_1,\ldots,\alpha_N$, respectively, such
that the twisted derived categories $\Db(X_i,\alpha_i)$,
$i=1,\ldots,n$, are all Fourier--Mukai equivalent.
\end{proposition}

\begin{proof}
Let $p_1,\ldots,p_N$ be the first $N$ primes and consider the
following diagonal positive definite $(2\times2)$-matrices:
$C_N:={\rm diag}(4\prod p_i^2 ,4\prod p_i^2)$ and $B_i:={\rm
diag}(4p_i^2,4p_i^2)$. We denote the lattices defined by them as
$T$ respectively $T_i$ and their natural generator by $e_1,e_2$
and $f_{i,1},f_{i,2}$, respectively. Clearly, the lattice $T$ can
be embedded (non-primitively) into each of the $T_i$ by
$e_j\mapsto (\prod_{k\ne i} p_k)f_{i,j}$.

Using primitive embeddings of $T$ and $T_i$ into $\Lambda$ and the
surjectivity of the period map, one finds K3 surfaces $Z$ and
$X_i$ realizing $T$ respectively $T_i$ as their transcendental
lattices. In particular, the Picard number of all of them is $20$.

The methods of Proposition \ref{twistedisuntwisted} apply and yield
equivalences $\Db(Z)\cong\Db(X_i,\alpha_i)$. As $|{\rm
disc}(T(X_i))|=16p_i^4$, all K3 surfaces $X_i$ are pairwise
non-isomorphic.
\end{proof}

The above proof can be modified for other Picard numbers ($\rho\geq
8$) and with additional geometric conditions on the surfaces $X_i$,
e.g.\ to be elliptic or Kummer. We leave this to the reader.

\begin{remark}\label{notthesameorder}
The examples constructed above are also interesting from another
point of view. Namely, they provide examples of derived equivalent
twisted K3 surfaces $(X_1,\alpha_1)$ and $(X_2,\alpha_2)$ with
${\rm ord}(\alpha_1)\ne{\rm ord}(\alpha_2)$. Compare the
discussion in Remark \ref{Brauerorderdifferent}.
\end{remark}
 
\end{document}